\documentclass[11pt]{amsart}
\usepackage{amssymb,amsmath}
\title{A berndtsson-Andersson operator solving $\overline{\partial}$-equation with $W^\alpha$-estimates on convex domains of finite type}
\author{William ALEXANDRE}
\address{Laboratoire Paul Painlev\'e U.M.R. CNRS 8524, U.F.R. de Math\'ematiques, Universit\'e Lille 1, F59 655 Villeneuve d'Ascq Cedex, France.} \email{ william.alexandre@math.univ-lille1.fr}
\date{}
\newtheorem{theorem}{Theorem}[section]
\newtheorem{lemma}[theorem]{Lemma}
\newtheorem{proposition}[theorem]{Proposition}
\newtheorem{corollary}[theorem]{Corollary}
\newtheorem{remarque}{\it Remark}
\newtheorem{definition}[theorem]{Definition}
\def \pint {\vbox{ \hbox to 5 pt {\hfil \vrule height 4pt}\hrule}\hskip 3pt}
\def\leqs{\lesssim}
\def\geqs{\gtrsim}
\def\eqs{\eqsim}
\def\cc{\mathbb{C}}

\def\nn{\mathbb{N}}

\def\re{{\rm Re} }
\def\im{{\rm Im} }
\def \qed {\hbox{\hskip 5pt} \vbox{\hrule \hbox to 5pt 
{\vrule height 4.2pt \hfil \vrule}\hrule}}
\def \pint {\vbox{ \hbox to 5 pt {\hfil \vrule height 4pt}\hrule}\hskip 3pt}

\newcommand{\cal}{\mathcal}
\newcommand{\wrt}{with respect to }

\newcommand{\mat}[3]{\left #1  #2  \right #3}
\newcommand{\diffp}[2]{\frac{\partial #1}{\partial #2}}
\newcommand{\mlabel}[1]{\label {#1}}
\renewcommand{\over}[2]{\genfrac{}{}{0pt}{}{#1}{#2}}
\newcommand{\dbar}{\overline\partial}
\newcommand{\w}[2]{W^{#1}_{0,#2}}
\newcommand{\p}[2]{{\cal P}_{#1}(#2)}
\newcommand{\po}[2]{{\cal P}^{(0)}_{#1}(#2)}

\newcommand{\pr} {{\it Proof:} }

\begin{document}
\pagestyle{plain}

\begin{abstract}
The Carleson measures were first introduced by Carleson in order to solve the corona problem for the disc in $\cc$. The notion of Carleson measure can be generalized to any homogeneous space and were also used in the context of the corona problem in $\cc^n$ for example by Varopoulos,  Amar and Andersson and Carlsson. One of the steps to solve the $H^p$ corona problem in a pseudoconvex domain is to  solve the $\overline\partial$ equation for a form $\mu$ satisifying a Carleson condition and get norm estimates of the solution in term of the Carleson norm of $\mu$. The main goal of this paper is to consider this question in the case of convex domains of finite type and to get estimates linked to the multitype of the domain.
\end{abstract}
\maketitle
\section{Introduction}
Our first interest for Carleson measures and $W^\alpha$-estimates on convex domains of finite type comes from the corona problem. This problem can be formulated for any domain $D$ of $\cc^n$, $n\geq 1$ as follows:
Let $g_1,\ldots, g_k$ be bounded holomorphic functions on $D$ such that $\sum_{j=1}^k \left| g_j \right|^2 \geq \delta>0$ on $D$. Given $\phi$ bounded and holomorphic on $D$, do there exist $f_1,\ldots, f_k$ bounded and holomorphic on $D$ such that $\sum_{j=1}^k f_jg_j=\phi$?

This problem was solved for the unit disc in the complex plan in \cite{Car} where Carleson introduced for the first time the notion of Carleson measure. However the corona problem is still unsolved for $n>1$. It is known that the answer is negative for some pseudoconvex domains of $\cc^2$ (see \cite{FSi,Sib}) but there are partial results in the positive direction. For example if $D=\{z\in\cc^n, r(z)<0\}$ let $H^p(D)$, $1\leq p<\infty$ be the set of all holomorphic functions $f$ over $D$ such that $\|f\|_p:= \left( \sup_{\varepsilon>0} \int_{bD_{-\varepsilon}} |f|^p d\sigma \right)^{\frac{1}{p}}<\infty$ where $D_{-\varepsilon}=\{z\in\cc^n, r(z)<-\varepsilon\}$, $bD_{-\varepsilon}$ is the boundary of $D_{-\varepsilon}$ and $\sigma$ is the area measure on $bD_{-\varepsilon}$. Amar proved in \cite{Ama} when $k=2$ and when $D$ is the ball of $\cc^n$ that if $\phi\in H^p(D)$ then there exist $f_1,f_2\in H^p(D)$ such that $f_1g_1+f_2g_2=\phi$. This result was generalized by Lin in \cite{Lin} for any $k$ in the polydisc of $\cc^n$ and by Andersson and Carlsson in \cite{AC0,AC0bis, AC} for any $k$ and for $D$ a general strictly pseudoconvex domain of $\cc^n$.

One strategy to solve the $H^p$-corona problem uses the Koszul complex: Let $D$ be a bounded domain of $\cc^n$ with smooth boundary, let $\Lambda ^l$ denotes the set of elements $f$ of degree $l$ of the exterior algebra of basis $e_1,\ldots, e_k$: $f=\sum_{i_1<\ldots<i_l} f_{i_1,\ldots, i_l} e_{i_1}\cap \ldots \cap e_{i_l}$ and let $e^*_1,\ldots, e^*_k$ be the dual basis of $e_1,\ldots, e_k$. We then define the mapping $\delta_g: H^p(\Lambda )\to H^p(\Lambda ^{l-1})$ by 
$\delta_g=\sum_{j=1}^k g_j e^*_j$. To solve the corona problem is equivalent to find $f=\sum_{j=1}^k f_je_j$ such that $\delta_gf=\phi $ with appropriate growth estimates for $f_1,\ldots, f_k$.\\
Let $\gamma_1,\ldots, \gamma_k$ be smooth functions over $D$ such that $\sum_{j=1}^k \gamma_jg_j=1$, that is $\delta_g(\gamma)=1$ where $\gamma=\sum_{j=1}^k \gamma_je_j$. We have the following Theorem (see \cite{AC}, Theorem 3.1).
\begin{theorem}
 Suppose that $T: C^\infty_*(\overline{D})\to C^\infty_*(\overline D)$ is such that $\overline\partial Tf=f$ if $\overline\partial f=0$, let $\phi $ be holomorphic with value in $\Lambda ^l$ such that $\delta_g(\phi )=0$ and let $r=\min(n,m-l-1)$. Then 
$f=\sum_{j=0}^r (-1)^j (\delta_g T)^j (\gamma \cap(\overline\partial \gamma)^j \cap \phi )$ is holomorphic and satisfies $\delta_gf=\phi $.
\end{theorem}
Hence such a function $f$ is a solution to the $H^p$-corona problem provided $f$ belongs to $H^p(D).$ This condition relies on the regularity of $T$. For the term 
$(\delta_g T) (\gamma \cap\overline\partial \gamma \cap \phi )$ one must show some W\"olf type estimates and for the term $(\delta_g T)^j (\gamma \cap(\overline\partial \gamma)^j \cap \phi )$, $j\geq 2,$ one must show estimates for Carleson measures. In this paper we are interested in these last estimates named $W^\alpha$-estimates.

Our operator $T$ which satisfies these estimates is a Berndtsson-Andersson integral operator. Such operators use weighted singular kernels whose singularity is given by a smooth function $H(\zeta,z)$ such that $H(\zeta,\zeta )=0$. The simplest choice is to set $H(\zeta ,z)=|\zeta -z|^2$ but one  can also use the more complicated choice $H(\zeta ,z)=|\zeta -z|^2 +|S(\zeta ,z)|^2$ where $S$ is a support function which reflects  the geometry of the domain (see \cite{ AC, BCD, Cum02, CF, DM, Ngu}). To fit the geometry of the boundary of the domain, one may think it would be more efficient to use only the support function. But this is impossible: one must add the part $|\zeta-z|^2$ in order to ensure the integrability of the kernel. And it seems that until now the term $|\zeta-z|^2$ was some kind of inert term ensuring the integrability and even sometimes disappearing on the boundary (see \cite{BCD, DM, Ngu}). However, in our case, such an operator constructed with $|\zeta-z|^2$ will not give the $W^\alpha$-estimates we are looking for because $|\zeta-z|^2$ is not linked to the more complicated geometry of the boundary. We will replace the term $|\zeta-z|^2$ by another one which is linked to the Bergman metric. This term will both ensure the integrability of the kernel and give the $W^\alpha$-estimates thanks to the links of the Bergman metric with the geometric properties of the boundary (see \cite{McN,McN1, McN2}).

An other problem will be encountered: In order to get the $W^\alpha$-estimates even for non smooth forms, that is for currents,  we will use a definition of Carleson currents related to smooth vector fields. And since the $\overline\partial$-operator is linked to the geometry of the boundary of the domain, we have to be able to find smooth vector fields which describe the geometry of the boundary. However the known tools for convex domains, that is for example the extremal basis, the Yu basis, which give a precise description of the boundary and should  intuitively  be the best candidates are not smooth ! (see \cite{Hef2}). Therefore we will have to find  good smooth vector fields. We will define them using again the Bergman metric.

This article is organized as follows: In Section 2 we fixe our notations and state our main results: the $W^\alpha$-estimates depending on the multitype of the domain (see Theorem \ref{th1} and \ref{th2}). In Section 3 we recall some properties of convex domains of finite type and all the tools needed to construct the operator $T$. In section \ref{section3} we give the estimates of all the parts of the kernel. In Section \ref{section4} we establish the $W^\alpha$-estimates and in section 6 we prove Theorem \ref{th2}.
\section{Notations  and main results}
In order to be as clear as possible we divide this section in the following subsection. In Subsection \ref{sub2.1} we define the multitype of a convex domain of finite type. In Subsection \ref{sub2.2} we define the homogeneous structure of the boundary of $D$ and the Carleson measures. In Subsection \ref{sub2.3} we state our main results.
\subsection{Multitype of a convex domain of finite type} \mlabel{sub2.1}
We consider a bounded convex domain $D:=\{z\in\cc^n,\ r(z)<0\}$ of finite type $m$ with $C^\infty$ smooth boundary, $r$ a smooth convex function whose gradiant does not vanish in a neighborhood of $bD$ the boundary of $D$.
Let $D_\alpha$ denote the set $\{z\in\cc^n,\ r(z)<\alpha\},$ let $\eta_z$ be the outer unit normal to $bD_{r(z)}$, the boundary of $D_{r(z)}$, at the point $z$ and let $T^\cc_z bD_{r(z)}$ denote the complex tangent space at $z$.\\
For $f:\cc\to \cc$ such that $f(0)=0$, we denote by $\nu (f)$ the multiplicity of $0$ as a zero of $f$.
\begin{definition}
The variety 1-type $\Delta_1(bD,p)$ of $bD$ at a point $p$ is defined as
$$\Delta_1(bD,p)=\sup_z\frac{\nu(z^*r)}{\nu(z-p)}$$
where the supremum is taken over all non zero germ $z:\Delta \to\cc^n$ from $\Delta $, the unit disc in $\cc$, into $\cc^n$, such that $z(0)=p$. The function $z^*r$ is the pullback of $r$ by $z$.\\
The variety $q$-type $\Delta _q(p,bD)$ at the point $p$ is defined as
$$\Delta _q(bD,p):=\inf_H \Delta _1(bD\cap H,p)$$
where the infimum is taken over all $(n-q+1)$-dimensional complex linear manifolds $H$ passing through $p$.
\end{definition}
In the case of smooth convex domain, one can define the multitype of $D$ in the following way (see \cite{BS,McN,Yu}).
\begin{definition}
 The multitype ${\rm M}(bD,p)$ of $bD$ at the point $p$ is defined to be the $n$-tuple $(\Delta _{n}(bD,p),\Delta _{n-1}(bD,p),\ldots, \Delta _{1}(bD,p))$.\\
 The multitype ${\rm M}(bD)$ of $bD$ is the $n$-tuple $${\rm M}(bD):=(\sup_{p\in bD}\Delta _{n}(bD,p),\ldots, \sup_{p\in bD}\Delta _{1}(bD,p)).$$
\end{definition}
The type of $D$ is in fact the last entry of $M(bD)$: $m=\sup_{p\in bD} \Delta_1(bD,p)$.
\subsection{Carleson measure on a convex domain of finite type} \mlabel{sub2.2}
The notion of Carleson measure can be defined on any space endowed with a structure of homogeneous space (see \cite{AB,CW}). For a convex domain of finite type, this structure is  induced by the polydiscs of McNeal. They are defined as follows (see \cite{BCD,McN, McN1}). For a point $z$ near $bD$ and for a sufficiently small $\varepsilon>0$ we set $$\tau(z,v,\varepsilon):=\sup\{t>0, |r(z+\lambda v)-r(z)|<\varepsilon, \forall \lambda \in\cc,\ |\lambda |<t\}.$$
In other words $\tau(z,v,\varepsilon )$ is the distance from $z$ to the level set $\{r=r(z)+\varepsilon\}$ in the complex direction $v$.\\
We denote by $w^*_1,\ldots, w^*_n$ an $\varepsilon$-extremal basis at $z$ as defined in \cite{BCD}. Such a basis is defined as follows:  $w_1^* =\eta_z$ and if $w^*_1,\ldots, w^*_{i-1}$ are already defined, then $w_i^*$ is a unit vector orthogonal to $w_1^*,\ldots, w^*_{i-1}$ such that $\tau(z,w^*_i,\varepsilon)= \sup_{\over{v\perp w^*_1,\ldots, w^*_{i-1}}{\|v\|=1}} \tau(z,v,\varepsilon)$.
We write $\tau_i(z,\varepsilon)=\tau(z,w_i^*,\varepsilon)$, for $i=1,\ldots,n$, and set $${\cal P}_\varepsilon(z):=\left\{ \zeta=z+\sum_{i=1}^n \zeta^*_iw_i^*\in\cc^n, |\zeta^*_i|<\tau_i(z,\varepsilon),\:i=1,\ldots n\right\}.$$ 
The two following propositions, proved in \cite{McN1}, show that the polydiscs define a structure of homogeneous space on $D$. 
\begin{proposition}\mlabel{prop1.1}
 For all $c>0$ there exists $b=b(c)>0$ such that for all $\varepsilon>0$ and all $z$ in a neighborhood of $bD$
\begin{align*}
{\cal P}_{c\varepsilon}(z)&\subset b{\cal P}_{\varepsilon}(z) \\
{c\cal P}_{\varepsilon}(z)&\subset {\cal P}_{b\varepsilon}(z) 
\end{align*}
\end{proposition}
In particular, there exist $c_1,c_2>0$ which do not depend on $z$ nor on $\varepsilon $ such that  $\text{ Vol}({\cal P}_\varepsilon(z))\leq c_1 \text{ Vol}({\cal P}_{c\varepsilon}(z))$ and $\text{ Vol}({\cal P}_{c\varepsilon}(z))\leq c_2 \text{ Vol}({\cal P}_{\varepsilon}(z))$. Since we will frequently formulate such inequalities depending on constants, we will write $A\leqs B$ if there exists a constant $c>0$ such that $A\leq cB$. Each time we will indicate the dependance of the constant. We will write $A\eqs B$ if $A\leqs B$ and $B\leqs A$ both hold.
\begin{proposition}\mlabel{prop1.2}
 There exists $C>0$ such that for all $\varepsilon>0$ and all $z,\zeta$ in a neighborhood of $bD$ the following holds true: if ${\cal P}_\varepsilon(z)\cap {\cal P}_\varepsilon(\zeta)\neq \emptyset$ we have ${\cal P}_\varepsilon(z)\subset C{\cal P}_\varepsilon(\zeta)$
\end{proposition}
In particular if ${\cal P}_\varepsilon(z)\cap {\cal P}_\varepsilon(\zeta)\neq \emptyset$ then $\text{ Vol}({\cal P}_\varepsilon(z))\eqs\text{ Vol}({\cal P}_\varepsilon(\zeta))$ uniformly \wrt $\zeta,z$ and $\varepsilon$.\\
We set for $\zeta,z$ near $bD$
$$\delta(z,\zeta):=\inf\{\varepsilon >0,\ \zeta\in {\cal P}_\varepsilon(z)\}.$$
Proposition \ref{prop1.1} and \ref{prop1.2} show that $\delta$ is a pseudodistance. Actually the structure of homogenuous space defined by $\delta$ and the polydiscs are the generalisation of the Koranyi distance and balls of strictly pseudoconvex domains. Therefore the following definition of Carleson measure is also a generalisation of Carleson measures of strictly pseudoconvex domains.
\begin{definition}
We say that a positive finite measure $\mu$ on $D$ is a {\it Carleson measure} and we write $\mu\in W^1(D)$ if 
$$\|\mu\|_{W^1}:= \sup_{{\zeta\in bD},{\varepsilon>0}}\frac {\mu(\p\varepsilon\zeta\cap D)}{\sigma({\cal P}_\varepsilon(\zeta)\cap bD)}<\infty$$
where $\sigma $ denotes the aera measure on $bD$.
\end{definition}
\subsection{Main results} \mlabel{sub2.3}
In \cite{AC}, the authors define a norm for currents which takes into account the fact that the $\dbar$-operator behaves differently in the tangential and normal directions. When $\mu$ is a $(0,q)$-form with measure coefficients, the norm $\|\mu\|$ of $\mu$ in \cite{AC} satisfies $\|\mu\|\eqs |r|^{\frac12} |\mu|+ |\overline\partial r\wedge \mu|$ where $|\mu|$ is the absolute value of $\mu$. This equality says that the tangential components of $\mu$, that is $\overline\partial r\wedge \mu$, are requiered to be less regular than the normal component. The difference of regularity is given by $|r|^{\frac12}$ where $\frac12={1-\frac12}$, that is 1 minus the order of contact of a tangent vector $v$ and the boundary of the domain. It is well known that the $\dbar$-operator behaves differently in a direction accordingly to the order of contact of that direction and of the boundary of the domain. So we will use a metric which take that fact into account. For $\varepsilon >0$, $z\in\cc^n$ and $v$ a non zero vector we set
\begin{align*}
k(z,v)&:=\frac{d(z)}{\tau(z,v,d(z))}
\end{align*}
where $d(z)=|r(z)|$.
For a fixed $z$, the convexity of $D$ implies that the function defined by $v\mapsto k(z,v)$ if $v \neq 0$, $0$ otherwise is a kind of non-isotropic norm. In fact, as the vectorial norm used in \cite{AC}, $k(z,\cdot)$ is equal to $d(z)$ times the Bergman metric.
The norm $k$ was already used in \cite{BCD} to define at every point $z\in D$ a punctual norm for $C^\infty$-smooth differential forms. When $\omega$ is a smooth 1-form on $\overline D$, Bruna, Charpentier and Dupain define $\|\omega(z)\|_{k}$ as the smooth function of $z$ by $\|\omega(z)\|_{k}:=\sup_{u\neq 0} \frac{|\omega(z)(u)|}{k(z,u)}$ which is the norm of the linear form $\omega(z)$ with respect to the norm $k(z,\cdot)$. When $u$ is a tangent vector such that the order of contact of $bD$ and
the line spanned by $u$ passing through $z$ is $m'$, then $k(z,u)\eqs d(z)^{1-\frac1{m'}}$ and the norm $\|\cdot\|_k$ quantify the difference of regularity as the norm of Andersson and Carlsson. 
However, in this paper we are interested in estimates which also generalize the results of Amar and Bonami in \cite{AB} for forms with measure coefficients which may be non smooth. Therefore such a definition does not make sense in our case and is forbidden. This is why we define the following norm. When $\mu$ is a $(0,q)$-current with measure coefficients, we can apply $\mu$ to $q$ smooth vector fields and we obtain in this way a measure. We can then define what we call a Carleson current.
\begin{definition}
We say that a $(0,q)$-form $\mu$ with measure coefficients  is a $(0,q)$-Carleson current if
$$\|\mu\|_{W^1_{0,q}}:=\sup_{u_1,\ldots,u_q} \left\|\frac1{k(\cdot,u_1)\ldots k(\cdot,u_q)}\left|\mu(\cdot)[u_1,\ldots,u_q]\right|\right\|_{W^1}<\infty.$$
where the supremum is taken over all smooth vector fields $u_1,\ldots, u_q$ which never vanish and where $|\mu(\cdot)[u_1,\ldots,u_q]|$ is the absolute value of the measure $\mu(\cdot)[u_1,\ldots,u_q]$.\\
We denote by $W^1_{0,q}(D)$ the set of all $(0,q)$-Carleson currents.
\end{definition}
Therefore $\|\cdot\|_{W^1_{0,q}}$ is a norm on forms with measure coefficients associated to the vectorial norm $k$ and is in the same spirit than the norms used in \cite{AB} and \cite{AC} but $\|\cdot\|_{W^1_{0,q}}$ take into account the non isotropy of the boundary of the domain. Moreover, we should notice that our norm is weaker than the norm of Bruna, Charpentier and Dupain in the sense that $\|\mu\|_{\w11}\leqs \int_B\|\mu(\zeta )\|_kdV(\zeta)$ for all smooth $\mu$.\\
Let $W^0$ be the set of positive bounded measures on $D$. For $\mu\in W^0$, we put $\|\mu\|_{W^0}:=\mu(D)$. Analogously to $W^1_{0,q}(D)$ we defined $\w0q(D)$:
\begin{definition}
We say that $\mu$ is a $(0,q)$-current with bounded measure coefficients and we write $\mu\in W^0_{0,q}(D)$ if
$$\|\mu\|_{W^0_{0,q}}:=\sup_{u_1,\ldots,u_q} \left\|\frac1{k(\cdot,u_1)\ldots k(\cdot,u_q)}\left|\mu(\cdot)[u_1,\ldots,u_q]\right|\right\|_{W^0}<\infty,$$
where the supremum is taken over all smooth vector fields $u_1,\ldots, u_q$ which never vanish and where $|\mu(\cdot)[u_1,\ldots,u_q]|$ is the absolute value of the measure $\mu(\cdot)[u_1,\ldots,u_q]$.
\end{definition}
We point out that we can apply  the regularization argument of \cite{AC-1} with a non smooth current $\mu$. It gives a sequence of smooth currents $(\mu_k)_{k\in\nn}$ compactly supported in $\overline D$ and weakly converging to $\mu$ such that $\|\mu_k\|_{W^1}$ and $\|\mu_k\|_{W^0}$ are controlled respectively by $\|\mu\|_{W^1}$ and $\|\mu\|_{W^0}$. This regularization argument allows us to deduce results for non smooth current from the corresponding results for smooth compactly supported currents. It would not have been the case if we were using the norm $\|\cdot\|_k$. However, it induced a major difficulty when it comes to work with the norms $\|\cdot\|_{W^\alpha}$, $\alpha=0,1$. These norms are computed using smooth vector fields and in general the supremum should be achieved with extremal basis which are not smooth (see \cite{Hef2}). We shall use the Bergman metric to overcome this difficulty (see section \ref{section3}).\\
For all $\alpha\in ]0,1[$ the space $\w\alpha q(D)$ will denote the complex interpolate space between $\w0q(D)$ and $\w1q(D)$. One can ``understand'' these spaces by the work  of  Amar and Bonami who proved in \cite{AB} that a measure $\mu$ belongs to $W^\alpha(D)$, $\alpha\in ]0,1[$, if and only if there exists a Carleson measure $\mu_1$ and $f\in L^{\frac1{1-\alpha}}(bD, d\mu_1)$ such that $\mu=fd\mu_1$.
We can now state our main results.
\begin{theorem}\mlabel{th1}
Let $D$ be a bounded convex domain of finite type with $C^\infty$-smooth boundary, let $(m_1,\ldots, m_n)$ be the multitype of $D$ and $\gamma_0>-1$. Then there exists a linear operator $T:C^\infty_*(\overline D)\to C^\infty_*(\overline D)$ such that 
\begin{trivlist}
 \item[$(i)$] $\overline\partial T\mu=\mu$ for all $\overline\partial$-closed $(0,q)$-forms $\mu$,
\item[$(ii)$] for all $\alpha\in[0,1]$, all $\gamma \in]-1,\gamma_0]$, all $q=1,\ldots, n$ and  all $\mu\in C^\infty_{0,q}(\overline D)$ such that $d(\cdot)^{q-\sum_{i=n-q+2}^n\frac1{m_i}+\gamma }\mu$ belongs to ${W^\alpha_{0,q}}(D)$: $$\left\|d(\cdot)^{q-1-\sum_{i=n-q+2}^n\frac1{m_i}+\gamma }T\mu\right\|_{W^\alpha_{0,q-1}}\leqs \left\|d(\cdot)^{q-\sum_{i=n-q+2}^n\frac1{m_i}+\gamma }\mu\right\|_{W^\alpha_{0,q}}$$
 uniformly with respect to $\mu$.
\end{trivlist}
\end{theorem}
In the first part of Theorem 4.1 of \cite{AC}, Andersson and Carlsson prove the same result for strictly pseudoconvex domains. The gain of regularity in our result and in the Andersson-Carlsson's result  are equal modulo the change of norm and the fact that the multitype of a strictly pseudoconvex domain is $(1,2,\ldots, 2)$. As in Theorem 4.1 of \cite{AC}, $\gamma =-1$ appears as a limit case but only when $q=1$. Before we state our result in this limit case, we recall the definition of BMO functions in convex domains of finite type.
\begin{definition}
 We say that $f\in L^1(bD)$ has Bounded Mean Oscillation and we write $f\in BMO(bD)$ when
$$\|f\|_*:=\sup_{\zeta\in bD,\varepsilon>0}\frac1{\sigma({\cal P}_\varepsilon(\zeta)\cap bD)}\int_{\p\varepsilon \zeta \cap bD} \left| f(\xi)-f_{\zeta ,\varepsilon}\right| d\sigma (\xi)<\infty,$$
where $f_{\zeta,\varepsilon}=\frac1{\sigma({\cal P}_\varepsilon(\zeta)\cap bD)}\int_{\p\varepsilon \zeta \cap bD} \left| f(\xi)\right| d\sigma (\xi)$.
\end{definition}
\begin{theorem}\mlabel{th2}
The operator of Theorem \ref{th1} satisfies for all smooth $(0,1)$-forms  $\mu\in\w\alpha1(D)$,
\begin{trivlist}
\item[$(i)$] $\left\|T\mu\right\|_{L^{p}(bD)}\leqs \left\|\mu\right\|_{W^\alpha_{0,1}}$, where $\frac1p=1-\alpha$, $\alpha\in[0,1[$,
\item[$(ii)$] $\left\|T\mu\right\|_{*}\leqs \left\|\mu\right\|_{W^1_{0,1}}$ if $\alpha=1$,
\end{trivlist}
 uniformly with respect to $\mu$.
\end{theorem}
This result is again a generalisation of the second part of Theorem 4.1 of \cite{AC} when $q=1$. 
A similar statement to Theorem \ref{th2} is obtained by Nguyen in \cite{Ngu}
but involving another operator $K$ and for the stronger norm $\|\cdot\|_k$ of \cite{BCD} for smooth forms. This operator $K$ is also a Berndtsson-Andersson operator but constructed without using the Bergman metric and it does not satisify the estimates of theorem \ref{th1}.\\
The regularization argument of \cite{AC-1} gives immediatly the following Corollaries:
\begin{corollary}\mlabel{cor9}
Let $D$ be a bounded convex domain of finite type with $C^\infty$-smooth boundary and let $(m_1,\ldots, m_n)$ be the multitype of $D$.\\ For all $\alpha\in[0,1]$, all $\gamma >-1$ and all $\overline\partial$-closed $(0,q)$-current $\mu$ such that $d(\cdot)^{q-\sum_{i=n-q+2}^n\frac1{m_i}+\gamma }\mu$ belongs to ${W^\alpha_{0,q}(D)}$, there exists a $(0,q-1)$-current $f$ such that 
\begin{trivlist}
 \item[$(i)$] $\overline\partial f=\mu$,
\item[$(ii)$]  $d(\cdot)^{q-1-\sum_{i=n-q+2}^n\frac1{m_i}+\gamma }f$ belongs to $W^\alpha_{0,q-1}(D)$ and satifies
$$\left\|d(\cdot)^{q-1-\sum_{i=n-q+2}^n\frac1{m_i}+\gamma }f\right\|_{W^\alpha_{0,q-1}}\leqs \left\|d(\cdot)^{q-\sum_{i=n-q+2}^n\frac1{m_i}+\gamma }\mu\right\|_{W^\alpha_{0,q}}$$
 uni\-form\-ly with respect to $\mu$.
\end{trivlist}
\end{corollary}
\begin{corollary}\mlabel{cor10}
Let $D$ be a bounded convex domain of finite type with $C^\infty$-smooth boundary.\\ For all $\alpha\in[0,1]$ and all $\overline\partial$-closed $(0,1)$-current $\mu\in {W^\alpha_{0,1}(D)}$ there exists a function $f$ such that 
\begin{trivlist}
\item[$(i)$] $\overline\partial f=\mu$ on $D$,
\item[$(ii)$] $\left\|f\right\|_{L^{p}(bD)}\leqs \left\|\mu\right\|_{W^\alpha_{0,1}}$, where $\frac1p=1-\alpha$, $\alpha\in[0,1[$,
\item[$(iii)$] $\left\|f\right\|_{*}\leqs \left\|\mu\right\|_{W^1_{0,1}}$ if $\alpha=1$,
\end{trivlist}
 uniformly with respect to $\mu$.
\end{corollary}
Corollary \ref{cor10} is a generalization of Theorem 7 in \cite{AB} in the case of $W^\alpha$ spaces.\\
We have chosen to formulate Theorem \ref{th1} and \ref{th2} and Corollary \ref{cor9} and \ref{cor10} with the norm $\|\cdot\|_{W^\alpha}$ associated to the vectorial norm $\|\cdot\|_k$. A more natural or more intrinsic choice could have been to use the Bergman metric instead of $\|\cdot\|_k$. The estimates would have been the same except that a factor $d(\cdot)$ to the power the degree of the form would have disappear.

\section{Construction of the operator}\mlabel{section2}
We first recall some properties of convex domains of finite type. 
We have (see \cite{Hef2}, Corollary 2.18)
\begin{proposition}\mlabel{prop3}
Let $z\in D$ be a point near $bD$, $(m_1,\ldots,m_n)$ denotes the multitype of $bD_{r(z)}$ at $z$, $\varepsilon>0$ and $w_1,\ldots, w_n$ is an $\varepsilon $-extremal basis at the point $z$. Then we have for $i=2,\ldots, n$ uniformly with respect to $z$ 
$$\tau_{i}(z,\varepsilon )\eqs \varepsilon ^{\frac1{m_{n-i+2}}}$$
and
$$\tau_1(z,\varepsilon)\eqs \varepsilon.$$
\end{proposition}
Among the extremal basis there exists a basis $w'_1,\ldots, w'_n$ of $\cc^n$ such that the order of contact of $bD$ and the line spanned by $w'_i$ passing through $z$ is equal to  $\Delta_{n+1-i}(bD,z)$. Such a basis is called a Yu basis at $z$ (see \cite{Hef}). This basis has the following properties.
\begin{proposition}\mlabel{prop4}
 Let $z\in D$ be a point near $bD$, $(m_1,\ldots,m_n)$ denotes the multitype of $bD_{r(z)}$ at $z$, $w'_1,\ldots, w'_n$ a Yu basis at $z$, $\varepsilon>0$ and $w^*_1,\ldots, w^*_n$ a $\varepsilon$-extremal basis at $z$ and $v= \sum_{j=1}^n v^*_j w^*_j=\sum_{j=1}^n v'_j w'_j$ a unit vector. Then, uniformly \wrt $z, v$ and $\varepsilon $ we have
\begin{align*}
\frac1{\tau (z,v,\varepsilon )}&\eqs \sum_{j=1}^n\frac{|v'_j|}{\varepsilon ^{\frac1{m_j}}}\eqs \sum_{j=1}^n\frac{|v^*_j|}{\tau_j(z,\varepsilon )}.
\end{align*}
\end{proposition}
The first ``equality'' is shown in \cite{Hef2}, Theorem 2.22, while the second is shown in \cite{McN}, Proposition 2.2. We notice that in particular, with the notations of Proposition \ref{prop4}, $\tau(z,w'_j,\varepsilon)\eqs \varepsilon^{\frac1{m_j}}$.\\
The two next properties are proved in \cite{BCD} and \cite{McN1} respectively.
\begin{proposition}\mlabel{prop1}
 Let $z\in D$ be a point near $bD$, $v$ a unit vector in $\cc^n$ and $\varepsilon _1\geq \varepsilon _2>0$. Then we have uniformly with respect to $z$, $\varepsilon _1,\varepsilon _2$ and $v$
$$\left(\frac{\varepsilon_1}{ \varepsilon _2}\right)^{\frac1m}\leqs\frac{\tau (z,v,\varepsilon _1)}{\tau (z,v,\varepsilon _2)}\leqs \frac{\varepsilon _1}{\varepsilon _2}.$$
\end{proposition}
%\begin{corollary}
%Let $z\in D$ be a point near $bD$ and  $K,\varepsilon>0$. Then, uniformly \wrt $\varepsilon$ but not to $K$ we have 
%$$Vol\bigl(\p\varepsilon z\bigl)\eqs Vol\bigl(\p{K\varepsilon} z\bigl)$$
%\end{corollary}
\begin{proposition}\mlabel{prop5}
 Let $z\in D$ be a point near $bD$, $v$ a unit vector in $\cc^n$, $\varepsilon>0$ and $\zeta\in\p\varepsilon z$. Then we have uniformly with respect to $z$, $\zeta$, $\varepsilon $ and $v$
$$\tau(z,v,\varepsilon)\eqs\tau(\zeta,v,\varepsilon).$$
\end{proposition}
We now recall the definition and some properties of the Bergman metric that we will need (see \cite{Ran}). 
The orthogonal projection from $L^2(D)$ onto $L^2(D)\cap {\cal O}(D)$, where $\cal O(D)$ is the set of holomorphic function on $D$, is called the Bergman projection. We denote it by $\cal B$. There exists a unique integral kernel $B$  such that for all $f\in L^2(D)$
$${\cal B}f(z)=\int_DB(\zeta,z) f(\zeta )dV(\zeta).$$
The kernel $B(\zeta,z)$ is call the Bergman kernel. This kernel is holomorphic with respect to $z$, antiholomorphic with respect to $\zeta$ and satisfies $B(\zeta ,z)=\overline{B(z,\zeta )}$.\\
The Bergman metric $M_B(z,\cdot)$ for $z\in D$ is an hermitian metric defined by the matrix $(b_{i,j}(z))_{i,j=1,\ldots, n}$ where $b_{i,j}(z)=\diffp{^2}{z_i\partial\overline{z_j}} \ln B(z,z)$. This means that the Bergman norm of $v=\sum_{i=1}^n v_i e_i$, where $e_1,\ldots, e_n$ is the canonical basis of $\cc^n$, is given by $M_B(z,v)=\left(\sum_{i,j=1}^n b_{ij}(z)v_j\overline{v_i}\right)^{\frac12}$.\\
Using Theorem 3.4 and 5.2 of \cite{McN1} and Proposition \ref{prop4}, we easily get
\begin{theorem}\mlabel{th3}
For all $z\in D$ in a neighborhood of $bD$ we have 
$$B(z,z) \geqs \frac1{Vol \bigl(\p{d(z)}z\bigr)}.$$
Let $w$ be any orthonormal coordinate system centered at $z$ and let $v_j$
be the unit vector in the $w_j$-direction. Then we have uniformly with respect to $z$
$$\left|\diffp{^{|\alpha|+|\beta |} B }{w^\alpha \partial \overline{w}^\beta} (z,z)\right|\leqs \frac1{Vol\bigl(\p{d(z)}z\bigr) \prod_{j=1}^n \tau(z,v_j,d(z))^{\alpha_j+\beta_j} }.$$
\end{theorem}
Theorem \ref{th3} yields to the following Corollary
\begin{corollary}\mlabel{cor1}
Let $z\in D$ be a point near $bD$, let $w$ be any orthonormal coordinate system centered at $z$, let $v_j$
be the unit vector in the $w_j$-direction and let $(b_{ij}^w)_{i,j}$ be the Bergman matrix in the $w$-coordinates. Then we have uniformly with respect to $z$
$$\left|\diffp{^{|\alpha|+|\beta|} b^w_{ij}
}{w^\alpha \partial \overline{w}^\beta} (z)\right|\leqs \frac1{\prod_{j=1}^n \tau(z,v_j,d(z))^{\alpha_j+\beta_j} }.$$
\end{corollary}
J.D. McNeal proved in \cite{McN2} that the eigenvalues of the matrix $\bigl(b_{ij}(z)\bigr)_{i,j}$ are $\tau _1(z,d(z))^{-2},\ldots, \tau_n(z,d(z))^{-2}$. Therefore we have
\begin{proposition}\mlabel{prop7}
Let $z\in D$ be a point near $bD$. Then uniformly with respect to $z$ 
$$\det\bigl((b_{i,j}(z)\bigr)\eqs \frac1{Vol (\p{d(z)}z)}.$$
\end{proposition}
The following  proposition is proved in \cite{McN2}.
\begin{proposition}\mlabel{prop2}
Let $z\in D$ be a point near $bD$, $v$ a unit vector in $\cc^n$. Then uniformly with respect to $z$ and $v$
$$M_B(z,v)\eqs \frac1{\tau (z,v,d(z))}.$$
\end{proposition}
As we explained in the last part of the introduction, we want to replace the term $|\zeta -z|^2$ used to ensure the integrability of the kernel in \cite{AC, DM, Ngu} by a term depending on the Bergman metric and vanishing on the boundary as in \cite{AC, DM, Ngu}. We consider $$\|v\|_{B,z}:=d(z)^2 M_B(z,v).$$
When $K$ is a compact subset of $D$, there exists $c_K>0$ such that for all $z\in K$  and all $v\in\cc^n$: $\|v\|^2_{B,z}\geq c_K|v|^2$. We also have for all $z_0\in bD$, $\lim_{z\to z_0, z\in D} \|v\|^2_{B,z} =0$.\\
We need a Hefer section for $\|\zeta-z\|^2_{B,z}$. If we define for $j=1,\ldots, n$
\begin{align*}
b_j(\zeta,z)&=d(z)^4 \sum_{i=1}^n b_{i,j}(z)\overline{\zeta_i-z_i}
\end{align*}
and $b(\zeta,z)=\sum_{j=1}^nb_j(\zeta ,z) d\zeta _j$, we get $\langle b(\zeta ,z), \zeta -z\rangle=\|\zeta-z\|_{B,z}^2$ where 
for $\alpha=\sum_{j=1}^n \alpha_jd\zeta_j$ we set $\langle \alpha,\zeta -z\rangle:= \sum_{j=1}^n\alpha_j(\zeta _j-z_j)$.\\
The $(0,1)$-form $b$ is of class $C^1$ on $\overline D\times \overline D$ thanks to the weight $d(z)^4$. Morevover for all $K$ compact set of $D$  there exists $c_K>0$ such that $$|\langle b(\zeta ,z), \zeta -z\rangle|\geq c_K |\zeta -z|^2.$$
The second ingredient of our kernel will be the support function of Diederich and Forn\ae ss constructed in \cite{DF}. We recall the definition of this support function. We fix some $\zeta$ in a neighborhood $\cal V$ of $bD$. We choose 
an orthonormal basis $w'_1,\ldots, w'_n$ such that $w'_1=\eta_\zeta$ and set $r_\zeta(\omega)=r(\zeta+\omega_1 w'_1+\ldots+\omega_n w'_n)$ and 
\begin{eqnarray*}
F_\zeta(\omega):=3\omega_1 + K\omega_1^2-K'\sum_{j=2}^m\kappa_j M^{2^j}\sum_{\over{|\beta|=j} {\beta_1=0}}\frac1{\beta!}\diffp{^jr_\zeta}{\omega^\beta}(0)\omega^\beta
\end{eqnarray*}
where $K,\:K',\:M$ are positive real numbers, $\kappa_j=1$ when $j\equiv\:0\: mod\:4$, $\kappa_j=-1$ when $j\equiv 2\:mod\:4$ and $\kappa_j=0$ otherwise.\\ 
We write $z\in \cc^n$ as $z=\zeta+\omega_{1,z} w'_1+\ldots+ \omega_{n,z} w'_n$ and define $F(\zeta,z)$ by 
$$F(\zeta,z):=F_\zeta(\omega_{1,z},\ldots, \omega_{n,z}).$$
$F$ satisfies the following theorem:
\begin{theorem} \mlabel{th03}
There exist a neighborhood $\cal V$ of $bD$ and positive constants $M$, $K$, $K',$ $k'$, $c_+$, $c_-$ and $R$  such that for all $\zeta\in{\cal V}$, all unit vector $v\in T^\cc_\zeta bD_{r(\zeta)}$ and all $w=(w_1,w_2)\in\cc^2$, with $|w|<R,$ we have 
\begin{eqnarray*}
\lefteqn{\re F(\zeta, \zeta+w_1\eta_\zeta+w_2v)\leq}\\
&\leq&-\mat{|}{\frac{\re w_1}2}{|}-\frac K2(\im w_1)^2-\frac{K'k'}4\sum_{j=2}^m\sum_{\alpha+\beta=j}\
\mat{|}{\mat{.} {\diffp{^jr(\zeta+\lambda v)}{\lambda^\alpha\partial \overline\lambda^\beta}}{|}_{\lambda=0}}{|}|w_2|^j\\
& &-c_\pm(r(\zeta)-r(\zeta+w_1\eta_\zeta+w_2v))
\end{eqnarray*}
where $c_{\pm}= c_-$ when $r(\zeta)-r(\zeta+w_1\eta_\zeta+w_2v)>0$ and $c_\pm=c_+$ otherwise.
\end{theorem}
This Theorem was proved in \cite{DF}. Dividing $F$ by $\frac1{2c_-}$ we can assume that $c_-=\frac12$ and $c_+\leq \frac12$. One should notice that we may have $F(\zeta,z)= 0$ when $|\zeta-z|>R$ so we must use a global version of this support function. For example we can construct such a function $S$ as in \cite{WA0}. This construction does not require other ideas than those of \cite{Ran}. As in the strictly pseudoconvex case (see \cite{Ran}) $S$ satisfies 
\begin{trivlist}
\item[i)] $S$ is of regularity $C^\infty$ in ${\cal V \times U}$, where $\cal U$ is a neighborhood of $\overline D$ and $S (\zeta, \cdot)$ is holomorphic on $\cal U$.
\item[ii)]  $S(\zeta,\zeta)=0$ for $\zeta\in {\cal U\cap V}$.
\item[iii)] There exists a constant $c>0$ such that $\re S (\zeta,z)\leq -c|\zeta-z|^m$ for all $(\zeta,z)\in \cal V \times U$ with $r(\zeta)\geq r(z)$.
\item[(iv)] 
On $\{(\zeta,z)\in {\cal V\times U}, |\zeta-z|<\frac R2\}$, there exists a $C^\infty$-function $A$ with $\frac12\leq |A(\zeta,z)|\leq\frac32$ and $S=A \cdot F$.
Moreover $A(\zeta,z)=\frac1{1+(m'-v(\zeta,z))F(\zeta,z)}$ where $m'$ is a constant and $v$ a bounded $C^\infty$ function defined on ${\cal V\times U}$ such that all its derivatives are also bounded on ${\cal V\times U}$.
\end{trivlist}
Define $Q_j(\zeta ,z)=-\int_0^1 \diffp{S}{z_j}(\zeta,\zeta+t(z-\zeta))dt$, $j=1,\ldots, n$ to be the Hefer decomposition of $S$ given in \cite{WA}. Therefore $Q_1,\ldots, Q_n$ satisfy $S(\zeta ,z)=\sum_{j=1}^n Q_j(\zeta ,z)(\zeta _j-z_j)$.

We now have all the tools we need to define our operator. Since $S$ is only defined for $\zeta $ in a neighborhood of $bD$, we have to  truncate $S$. Let $\eta_0>0$ be small enough so that $D_{\eta_0}\setminus D_{-\eta_0}$ is contained in the neighborhood $\cal V$ of Theorem \ref{th03}. Let $\chi$ be a smooth function with compact support such that $\chi\equiv 1$ on $D_{-\eta_0}$ and $\chi\equiv 0$ on $\cc^n\setminus D_{-\frac{\eta_0}2}$.
We set
\begin{align*}
 Q(\zeta ,z)=&\sum_{j=1}^n Q_j(\zeta ,z)d\zeta _j,\\
q(\zeta ,z)=&(1-\chi(\zeta )) Q(\zeta ,z)+\chi (\zeta )\partial r(\zeta ),\\
\tilde Q(\zeta ,z)=&\frac1{r(\zeta )} q(\zeta ,z),\\
\tilde S(\zeta ,z)=&\langle \tilde Q(\zeta ,z), \zeta _j-z_j\rangle,\\
s(\zeta ,z)=&|\langle q(\zeta ,z), \zeta-z\rangle|^2 \overline{\langle q(\zeta ,z), \zeta-z\rangle} q(\zeta ,z)\\
&+|\langle q(z,\zeta), \zeta-z\rangle|^2\overline{\langle q(z,\zeta), \zeta-z\rangle} q(z,\zeta),\\
\tilde s(\zeta ,z)=& b(\zeta,z)+ s(\zeta ,z).
\end{align*}
With such a choice of truncature $\tilde S(\zeta ,z)$ is essentially equal to  $\frac1{r(\zeta)} S(\zeta ,z)$ when $\zeta $ is near the boundary and $|\langle s(\zeta,z),\zeta-z\rangle|$ is essentially $\delta(\zeta,z)^4$ when $\zeta$ is close to $bD$ and $z$ close to $\zeta$. This will cause no problem because the properties of the $\overline \partial$ solving kernel are only important near the boundary. We now define our operator. Let
$$K(\zeta ,z)=\sum_{k=0}^{n-1} c_{n,k} \left(\frac1{1+\tilde S(\zeta ,z)}\right)^{N+k} \frac{\tilde s \wedge (\overline\partial_z \tilde s)^{q-1}\wedge (\overline \partial_\zeta  \tilde s)^{n-k-q}\wedge (\overline\partial_\zeta  \tilde Q)^k}{\langle \tilde s(\zeta ,z),\zeta -z\rangle^{n-k}},$$
where $c_{n,k}$ is a constant and $N\in\nn$ large enough.\\ The form $\tilde s$ is of class $C^1$ over $\overline D\times \overline{D}$ and for all $K$ compact subset of $D$ there exists $c_K>0$ such that $|\langle\tilde s(\zeta,z)\zeta-z\rangle|\geq c_K|\zeta-z|^2$ for all $z\in K$ and all $\zeta\in D$. Morevover $\tilde Q(\zeta,\cdot)$ is holomorphic when $\zeta$ is fixed in $D$. Hence the kernel $K$ satisfies the hypothesis of \cite{BA} and the operator $T$ defined for any $f\in C^1_{0,q}(\overline{D})$, $q=1,\ldots, n$, by
$$Tf(z):=\int_Df(\zeta )\wedge K(\zeta ,z)$$
 is such that $\overline \partial Tf=f$ for all $\overline\partial$-closed $f$. Therefore, in order to prove Theorems \ref{th1} and \ref{th2} it remains to prove that $Tf$ satisfies the estimates announced in these theorems.
\section{Estimates of the kernel}\mlabel{section3}
We now give the estimates of the different terms of the kernel that we will need. 
In order to get the best estimates for $K$ we fix $\zeta\in D$, choose a Yu basis $w'_1,\ldots, w'_n$ at $\zeta$ and write the kernel in this basis. We write $z'=(z'_1,\ldots, z'_n)$ the coordinates of a point $z$ in the coordinates system centered at $\zeta$ with respect to the basis $w'_1,\ldots, w'_n$ and set $\tau(\zeta,w'_i,\varepsilon)=\tau'_i(\zeta,\varepsilon)$. We write $b(\zeta,z)=\sum_{j=1}^n b'_j(\zeta,z)d\zeta'_j$, $s(\zeta,z)=\sum_{j=1}^ns'_j(\zeta,z)d\zeta'_j$, $\tilde s(\zeta,z)=\sum_{j=1}^n\tilde s'_j(\zeta,z)d\zeta'_j$ and $\tilde Q(\zeta,z)=\sum_{j=1}^n\tilde Q'_j(\zeta,z)d\zeta'_j$. We estimate these forms and their derivatives. 
We begin by recalling the following proposition (see \cite{WA,BCD,DFF}):
\begin{proposition}\mlabel{prop30}
 Let $w$ be any orthonormal coordinate system centered at $\zeta$ and let $v_j$
be the unit vector in the $w_j$-direction. For all multiindices $\alpha$ and $\beta $ with $|\alpha+\beta |\geq 1$ and all $z\in \p\varepsilon\zeta$:
$$\left|\diffp{^{|\alpha|+|\beta|} r}{w^\alpha\partial \overline w^\beta}(z) \right|\leqs \frac{\varepsilon }{\prod_{j=1}^n \tau (\zeta,v_j,\varepsilon)^{\alpha_j+\beta _j}}
$$
uniformly with respect to $z$, $\zeta$ and $\varepsilon$.
\end{proposition}
We denote by $k$ a positive constant such that 
for all $\zeta\in D$ sufficiently near the boundary, $\p{kd(\zeta)}\zeta\subset D_{-\frac{r(\zeta)}2}$ (see \cite{McN}).
\begin{proposition}\mlabel{prop9}
For all $\zeta\in  D$,  all $\varepsilon\geq k d(\zeta)$ and all $z\in\p\varepsilon\zeta$ we have uniformly
\begin{align*}
|b'_j(\zeta,z)|&\leqs \frac{\varepsilon^4}{\tau'_j(\zeta,\varepsilon)},\\
\left| \diffp{b'_j}{\overline{z'_k}}(\zeta,z) \right|,\left| \diffp{b'_j}{\overline{\zeta'_k}}(\zeta,z) \right|&\leqs \frac{\varepsilon^4}{\tau'_k(\zeta,\varepsilon)\tau'_j(\zeta,\varepsilon)}.
\end{align*}
Moreover if $z$ belongs to $\p{kd(\zeta)}\zeta$ we have uniformly
\begin{align*}
|b'_j(\zeta,z)|&\leqs  \frac{d(\zeta)^4}{\tau'_j(\zeta,d(\zeta))}\sum_{i=1}^n 
\frac{\left| z'_i \right|}{\tau'_i(\zeta,d(\zeta))} .
\end{align*}
\end{proposition}
\pr
In order to prove the first and the last inequalities at the same time we fix $\varepsilon \geq kd(\zeta )$. We denote by $\bigl( b'_{i,j}(z)\bigr)_{i,j}$ the matrix of the Bergman metric in the Yu basis $w'_1,\ldots, w'_n$. Thus for all $\xi\in D$ and all $z\in \p\varepsilon\zeta$ we have 
\begin{align*}
b'_j(\xi,z)&=d(z)^4 \sum_{i=1}^n b'_{i,j}(z) \overline{(\xi'_i-z'_i)}.
\end{align*}
Corollary \ref{cor1} gives $|b'_{ij}(z)|\leqs\frac{1}{\tau (z,w'_i,d(z))\tau (z,w'_j,d(z))}$. Since $z$ belongs to $\p\varepsilon\zeta$ we have $d(z)\leqs d(\zeta)+\varepsilon\leqs \varepsilon $. Proposition \ref{prop1} implies that $\frac{1}{\tau (z,w'_k,d(z))}\leqs \frac{\varepsilon }{d(z)}\frac{1}{\tau (z,w'_k,\varepsilon )}$ which together with Proposition \ref{prop5} yields to 
$\frac{1}{\tau (z,w'_k,d(z))}\leqs \frac{\varepsilon }{d(z)}\frac{1}{\tau'_k (\zeta,\varepsilon )}$ and thus 
$|b'_{ij}(z)|\leqs\left(\frac{\varepsilon }{d(z)}\right)^2 \frac{1}{\tau'_i(\zeta ,\varepsilon )\tau'_j(\zeta ,\varepsilon )}.$
Hence
\begin{align}
\left|d(z)^4 \sum_{i=1}^n b'_{i,j}(z) (\xi'_i-z'_i) \right| &\leqs \frac{\varepsilon ^4}{\tau '_j(\zeta ,\varepsilon )} \sum_{i=1}^n\frac{\left| \xi'_i-z'_i\right|}{\tau '_i(\zeta ,\varepsilon )},\mlabel{eq62}
\end{align}
When $\varepsilon =kd(\zeta )$ and $\xi=\zeta$ this proves the last inequality.\\
When $\varepsilon \geq kd(\zeta )$ and $\xi=\zeta$, since $z$ belongs to $\p\varepsilon \zeta $ we have $|\xi'_k-z'_k|=|z'_k|\leqs \tau '_k(\zeta ,\varepsilon )$ for $k=1,\ldots, n$ thus (\ref{eq62}) also implies the first inequality. The other inequalities of the proposition can be shown analogously.\qed
\begin{lemma}\mlabel{lem3}
For all $\zeta\in D$ close enough to $bD$, all $z\in \p{kd(\zeta)}\zeta$ we have
$$\left| \langle \tilde s(\zeta,z),\zeta-z\rangle \right|^{\frac{1}{2}}\geqs d(\zeta)^2\sum_{j=1}^n \frac{\left| z'_j \right|}{\tau'_j(\zeta,d(\zeta))}.$$
\end{lemma}
\pr If $\zeta$ is close enough to $bD$: $\left| \langle \tilde s(\zeta,z),\zeta-z\rangle \right|^{\frac{1}{2}}\geqs d(z)^2 M_B(z,\zeta-z)$.
Proposition \ref{prop2} leads to $d(z)^2 M_B(z,\zeta-z) \eqs \frac{d(z)^2} {\tau(z,\zeta-z,d(z))}\eqs \frac{d(z)^2|\zeta-z|}{\tau\left(z,\frac{\zeta-z}{\left| \zeta-z \right|},d(z)\right)}$.\\
Since $z$ belongs to $\p{kd(\zeta)}\zeta$ we have $d(z)\eqs d(\zeta)$ and by Propositions \ref{prop1} and \ref{prop4} we have $\frac{|\zeta-z|}{\tau\left(z,\frac{\zeta-z}{\left| \zeta-z \right|},d(z)\right)} \eqs \frac{|\zeta-z|}{\tau\left(\zeta,\frac{\zeta-z}{\left| \zeta-z \right|},d(\zeta)\right)}\eqs \sum_{j=1}^n \frac{\left| z'_j \right|}{\tau'_j(\zeta,d(\zeta))}$ which finishes the proof of the lemma.\qed

We also need estimates for the support function $S$. As in \cite{WA2,DFF, Hef,Hef2,Ngu} we shall use a covering of the domain of integration by poly annuli of the form $\po\varepsilon \zeta:= \p\varepsilon\zeta\setminus c\p\varepsilon\zeta$ where $c>0$ is such that $c\p\varepsilon\zeta$ is included in $\p{\frac12\varepsilon}\zeta$ for all $\zeta$ and all $\varepsilon>0$.
The required estimates for $S$ are given by the following.
\begin{proposition}\mlabel{prop11}
\begin{trivlist}
\item[$(i)$] For all $\zeta\in D$, all $z\in\overline {D_{r(\zeta)}} \cap \bigl(\p\varepsilon\zeta\setminus c\p\varepsilon\zeta\bigr)$ we have uniformly with respect to $\zeta,z$ and $\varepsilon$
$$|S(\zeta,z)|\eqs \varepsilon.$$
\item[$(ii)$] For all $\zeta\in D$, all sufficiently small $\varepsilon$ and all $z\in\overline D\cap \p\varepsilon\zeta$ we have uniformly with respect to $\zeta,z$ and $\varepsilon$
$$|r(\zeta)+S(\zeta,z)|\eqs \varepsilon+d(\zeta)+d(z).$$
\end{trivlist}
\end{proposition}
\pr The inequality $|S(\zeta,z)|\leqs \varepsilon$ is a consequence of Proposition \ref{prop30} and one should notice that there is no need to assume $r(z)\leq r(\zeta)$. The lower bound $\varepsilon\leqs \left| S(\zeta,z) \right|$ of $(i)$ was proved in \cite{WA}.\\
The inequality $|r(\zeta)+S(\zeta,z)|\leqs \varepsilon+d(\zeta)+d(z)$ also comes from the Proposition \ref{prop30}. To prove the lower bound of $(ii)$, one should notice that $S(\zeta,z)=\frac1{\left| 1+m'(1-v(\zeta,z))F(\zeta,z) \right|^2} \bigl(F(\zeta,z)+\overline{(m'-v(\zeta,z))}\left| F(\zeta,z) \right|^2\bigr)$ and since $ |F(\zeta,z)|\eqs |S(\zeta,z)|$ is controled by $\varepsilon$, it suffices to prove that for sufficiently small $\varepsilon>0$, $|r(\zeta)+F(\zeta,z)|\geqs \varepsilon+d(\zeta)+d(z).$\\
We have $|r(\zeta)+F(\zeta,z)|\geqs -r(\zeta)-\re F(\zeta,z)+\left| \im F(\zeta,z) \right|$ and  
according to Theorem \ref{th03}, if $r(\zeta)\geq r(z)$, when we write $z=\zeta+\lambda\eta_\zeta+\mu v$, $v\in T^\cc_\zeta bD_{r(\zeta)}$:
\begin{align*}
-r(\zeta)-\re F(\zeta,z)&\geqs \mat{|}{\frac{\re \lambda}2}{|}+(\im \lambda)^2+\sum_{j=2}^m\sum_{\alpha+\beta=j}\
\mat{|}{\mat{.} {\diffp{^jr(\zeta+\tilde\mu v)}{\tilde\mu^\alpha\partial \overline{\tilde\mu}^\beta}}{|}_{\tilde\mu=0}}{|}|\mu|^j\\
& +(c_+-1) r(\zeta)-c_+r(\zeta+w_1\eta_\zeta+w_2v))\\
&\geqs \mat{|}{\frac{\re \lambda}2}{|}+(\im \lambda)^2+\sum_{j=2}^m\sum_{\alpha+\beta=j}\
\mat{|}{\mat{.} {\diffp{^jr(\zeta+\tilde\mu v)}{\tilde\mu^\alpha\partial \overline{\tilde\mu}^\beta}}{|}_{\tilde\mu=0}}{|}|\mu|^j\\
&+d(\zeta)+d(z)
\end{align*}
The same inequality also hold when $r(\zeta)\leq r(z)$. Now, as shown in \cite{WA}, 
for sufficiently small $\varepsilon>0$ and sufficiently small $\tilde c>0$ if $|\mu|<\tilde c\tau(\zeta,v,\varepsilon)$ we have $|\lambda|\eqs \varepsilon$, uniformly with respect to $\zeta$, $z$, $\varepsilon$ and $\tilde c$.\\
If we assume $|\mu|< \tilde c\tau(\zeta,v,\varepsilon)$ then  $|\lambda|\eqs \varepsilon$ and thus Proposition 3.1 $(vi)$ of \cite{DFF} gives that
\begin{align*}
|r(\zeta)+F(\zeta,z)|&\geqs
  d(\zeta)+d(z) +\hskip -1pt |\lambda| -\hskip -1pt |\lambda|^2-\hskip -1pt  \sum_{j=2}^m
\sum_{\alpha+\beta=j} \mat|{ \mat.{\diffp{^jr(\zeta+\mu v)}{\tilde\mu^\alpha\partial \overline{\tilde\mu}^\beta}}|_{\tilde\mu=0} }||\mu|^j\\
&\geqs \varepsilon(1-\varepsilon+\tilde c)+d(\zeta)+d(z)
\end{align*}
and if $\tilde c$ and $\varepsilon$ are sufficiently small 
$|r(\zeta)+F(\zeta,z)|\geqs\varepsilon+d(\zeta)+d(z).$\\
Now, if $|\mu|\geq \tilde c\tau(\zeta,v,\varepsilon)$ by Proposition 3.1 $(vi)$ of \cite{DFF}
\begin{align*} 
|r(\zeta)+F(\zeta,z)|
&\geq\sum_{j=2}^m\sum_{\alpha+\beta=j}\
\mat{|}{\mat{.} {\diffp{^jr(\zeta+\tilde\mu v)}{\tilde\mu^\alpha\partial \overline{\tilde\mu}^\beta}}{|}_{\tilde\mu=0}}{|}|\mu|^j
+d(\zeta)+d(z)\\
&\geqs \varepsilon+d(\zeta)+d(z).\qed
\end{align*}
We will use the following Corollary
\begin{corollary}\mlabel{cor5}
For all $\zeta\in D$ sufficiently close to $bD$, all $\varepsilon>0$ and all $z\in \po\varepsilon\zeta$ we have uniformly
$$\left| \langle\tilde s(\zeta,z),\zeta-z\rangle \right|^{\frac{1}{2}}\geqs \varepsilon^2.$$
\end{corollary}
\pr we have $\left| \langle\tilde s(\zeta,z),\zeta-z\rangle \right|\geq \left|\langle q(\zeta,z),\zeta-z\rangle  \right|^4 + \left| \langle q(z,\zeta),z-\zeta\rangle \right|^4$ and according to Proposition \ref{prop11} we have $\left|\langle q(\zeta,z),\zeta-z\rangle  \right|\geqs \varepsilon$ if $r(\zeta)\geq r(z)$ and $\left|\langle q(z,\zeta),z-\zeta\rangle  \right|\geqs\varepsilon$ if $r(z)\geq r(\zeta)$.\qed
\begin{remarque}\mlabel{rem1}
Up to uniform constants, saying that $z$ belongs to $\po\varepsilon\zeta$ is equivalent to saying that $\delta(\zeta,z)\eqs \varepsilon$ so the Proposition \ref{prop11} and Corollary \ref{cor5} could have been formulated with $\delta(\zeta,z)$ instead of $\varepsilon$.
\end{remarque}
We also need estimates for the Hefer decomposition of $S$. As for $b$, we write $Q$ in the Yu basi at $\zeta$: $Q(\zeta,z)=\sum_{j=1}^nQ'_j(\zeta,z)d\zeta'_j$ and use the following proposition which comes from \cite{WA3,WA} and Proposition \ref{prop4}:
\begin{proposition}\mlabel{prop6}
For all $\zeta$ near enough $bD$, all sufficiently small $\varepsilon>0$, all $,\xi,z\in\p\varepsilon\zeta$ and $i, j,k =1,\ldots, n$ we have uniformly with respect to $\zeta$, $z$, $\xi$  and $\varepsilon$ 
\begin{align*}
|Q'_i(\xi,z)|&\leqs \frac{\varepsilon}{\tau'_j(\zeta,\varepsilon)},\\
\left| \diffp{Q'_i} {{z'_j}} (\xi,z
) \right|,\left| \diffp{Q'_i} {{\zeta'_j}} (\xi,z) \right|, \mat|{
\diffp{Q'_i} {\overline {\zeta'_j}} (\xi,z) }|&\leqs \frac{\varepsilon} {\tau'_j(\zeta,\varepsilon)\tau'_i(\zeta,\varepsilon)},\\
\left| \diffp{Q'_i}{\zeta'_j\partial z'_k}(\xi,z)\right|,\left| \diffp{Q'_i}{\overline{\zeta'_j}\partial z'_k}(\xi ,z)\right|&\leqs \frac{\varepsilon }{\tau'_i(\zeta ,\varepsilon )\tau'_j(\zeta ,\varepsilon )\tau'_k(\zeta ,\varepsilon )}.
\end{align*}
\end{proposition}
We deduce from Proposition \ref{prop4}, \ref{prop9}, \ref{prop11} and  \ref{prop6} the following corollary.
\begin{corollary}\mlabel{cor3}
 For all $\zeta\in D$ near enough $bD$, $\varepsilon\geq kd(\zeta )$ let $z$ be a point in $\p\varepsilon\zeta\cap D$, $w''_1,\ldots, w''_n$ be an orthonormal basis and denote by $(\xi''_1,\ldots, \xi''_n)$ the coordinates of $\xi\in \cc^n$ in the coordinate system centered in $z$ with respect to the  basis $w''_1,\ldots, w''_n$. Then
\begin{align*}
\left| {\tilde s'_i} (\zeta,z)\right|&\leqs \frac{\varepsilon^4}{\tau'_j(\zeta,\varepsilon)},\\
\left| \diffp{\tilde s'_i}{\overline{z''_j}} (\zeta,z)\right|&\leqs  \frac{\varepsilon^4}{\tau(\zeta,w''_j,\varepsilon)\tau'_j(\zeta,\varepsilon)},\\
\left| \diffp{\tilde s'_i}{\overline{\zeta'_j}} (\zeta,z)\right|&\leqs \frac{\varepsilon^4}{\tau'_j(\zeta,\varepsilon)\tau'_j(\zeta,\varepsilon)},\\
\left| \diffp{\tilde Q'_i}{\overline{\zeta'_j}} (\zeta,z)\right|&\leqs \left(\frac\varepsilon {d(\zeta )}\right)^2 \frac{\varepsilon}{\tau'_i(\zeta,\varepsilon)\tau'_j(\zeta,\varepsilon)}.\\
\end{align*}
For all $z\in\p{kd(\zeta)}\zeta$
$$\left| \tilde s'_i(\zeta,z) \right| \leqs \frac{d(\zeta)^4}{\tau'_i(\zeta,d(\zeta))} \sum_{j=1}^n \frac{\left| z'_j \right|}{\tau'_j(\zeta,d(\zeta))}.$$
\end{corollary}
On the one hand the hypothesis on the current $\mu$ in Theorems \ref{th1} and \ref{th2} are  related to the norm $\|\cdot \|_{\w\alpha q}$ which is a supremum over smooth vector fields. On the other hand the kernel $K$ is linked to geometry of the boundary of the domain. So we will need to find smooth vector fields linked to the geometry of the boundary in order to evaluate the exterior product $\mu\wedge K$. A natural choice would be to consider the Yu basis or the extremal basis as vector fields. These basis would be a good choice if they were smooth but they are not smooth in general as shown by the example of Hefer in \cite{Hef2}. To find good vector fields, we use the Bergman metric again.

Let us point out the following fact. Let $\|\cdot\|$ be an hermitian metric on $\cc^n$ defined by the positive definite hermitian matrix $b$. There exists a matrix $a$ such that $a^{-2}=b$. The linear map associated with the matrix $a$ sends the unit ball of $\cc^n$ equiped with the standard hermitian structure to the unit ball of the metric $\|\cdot\|$. We apply this fact to the Bergman metric $M_B(\zeta ,\cdot)$ for $\zeta\in D$. Proposition \ref{prop2} implies that the McNeal polydisc $\p{d(\zeta )}\zeta $ is almost the unit ball of the metric $M_B(\zeta ,\cdot)$. Therefore if $a(\zeta )$ satisfies $a(\zeta )^{-2} = \bigl(b_{i,j}(\zeta )\bigr)$, the column of $a(\zeta )$ are roughly speaking some kind of $d(\zeta )$-extremal basis at $\zeta $. This will give us our smooth vector field:

Let ${\cal H}_n^{++}$ denotes the manifold of positive definite hermitian matrices of size $n$ and let $\Phi: {\cal H}_n^{++}\to{\cal H}_n^{++} $ be defined for all $A\in{\cal H}^{++}_n$ by $\Phi(A)=A^{-2}$. $\Phi$ is a $C^\infty$-smooth diffeomorphism thanks to the inverse mapping theorem. Let $\varepsilon_1,\ldots, \varepsilon _n$ be the canonical basis of $\cc^n$ and set $e_k(\zeta )=\Phi^{-1}\bigl( (b_{i,j}(\zeta ))\bigr) \varepsilon _k$.\\
The basis $e_1(\zeta ),\ldots, e_n(\zeta )$ is an orthonormal basis for the Bergman metric and depends  smoothly on $\zeta$. These vector fields $e_1(\zeta ),\ldots, e_n(\zeta )$ are the smooth vector fields we are looking for. Let $(e'_{j1}(\zeta),\ldots, e'_{jn}(\zeta))$ denote the coordinates of $e_j(\zeta)$ in the Yu basis at $\zeta$. We estimate them in the following proposition:
\begin{proposition}\mlabel{prop10}
For all $\zeta\in D$ near enough the boundary of $D$, we have uniformly \wrt $\zeta$
\begin{align*}
|\det(e_1(\zeta),\ldots, e_n(\zeta))|&\eqs \left(Vol (\p{d(\zeta)}\zeta)\right)^{\frac12},\\
|e'_{ij} (\zeta)|&\leqs \tau'_j(\zeta,d(\zeta)).
\end{align*}
\end{proposition}
\pr The first estimate is a direct consequence of Proposition \ref{prop7}. In order to prove the second one we notice that $M_{B}(\zeta ,e_i(\zeta ))=1$ thus Proposition \ref{prop2} gives $1\eqs \sum_{j=1}^n \frac{|e'_{ij}(\zeta )|}{\tau'_j(\zeta ,d(\zeta ))}$ and so $|e'_{ij} (\zeta)|\leqs \tau'_j(\zeta,d(\zeta))$.\qed
\begin{remarque}
Despite the $\varepsilon $-extremal basis at $\zeta$ does not smoothly depend on $\varepsilon $ nor on $\zeta$, we can use the basis $e_1(\zeta ),\ldots, e_n(\zeta )$ in order to get a  description of $\p\varepsilon \zeta $ smoothly depending on $\zeta $ and $\varepsilon $. This can be done as follows. Without restriction we may assume that $0$ belongs to $D$. Let $p$ denotes the calibrator or gauge function for $D$, that is $p(\zeta )=\inf\{\lambda>0,\ \zeta \in \lambda D\}$ and take $r=p-1$ as a defining function for $D$. Since $p(\lambda \zeta)=\lambda p(\zeta )$, for all $\lambda>0$, all $\zeta $ close enough to $bD$, all unit vector $v$ and all $\mu\in\cc$, the point $\lambda \zeta +\lambda \mu v$ belongs to $\p{\lambda \varepsilon }{\lambda \zeta }$ if and only if $\zeta+\mu v$ belongs to $\p{\varepsilon }{\zeta }$. We set $\lambda_{\zeta ,\varepsilon }= (p(\zeta )-\varepsilon)^{-1}$ for $\zeta $ near $bD$. Then, since $e_1(\zeta), \ldots, e_n(\zeta )$ is an orthonormal basis for the Bergman metric, there exist $C>c>0$ not depending on $\zeta$ nor on $\varepsilon $ such that $c\p\varepsilon \zeta \subset \left\{\sum_{i=1}^n \alpha_i e_i\bigl(\lambda_{\zeta ,\varepsilon } \zeta \bigr),\ \sum_{i=1}^n |\alpha_i|^2\leq1\right\}\subset C\p\varepsilon \zeta $.
\end{remarque}
%We can now prove the $W^\alpha$-estimates.
\section{The $W^\alpha$-estimates}\mlabel{section4}
We first prove Theorem \ref{th1} for $\alpha=0$. Let $(m_1,\ldots, m_n)$ be the multitype of $D$, $\gamma\in]-1,\gamma_0]$, $f\in C^\infty_{0,q}(\overline{D})$ such that 
$d(\cdot)^{q-\sum_{i=n-q+2}^n\frac1{m_i}+\gamma }f$ belongs to ${W^0_{0,q}}(D)$ and let $u_1,\ldots, u_{q-1}$ be $q-1$ smooth never vanishing vector fields.\\
We show that
$\left\|d(\cdot)^{q-1-\sum_{i=n-q+2}^n\frac1{m_i}+\gamma}Tf(\cdot)\hskip-1pt \left[\frac{u_1(\cdot)}{k(\cdot, u_1(\cdot))},\ldots,\frac{u_q(\cdot)}{k(\cdot, u_q(\cdot))}\right]\right\|_{W^0}$
 is
do\-minated by $ \left\|d(\cdot)^{q-\sum_{i=n-q+2}^n\frac1{m_i}+\gamma }f\right\|_{W^0_{0,q}}$.
Since $f$ is a $(0,q)$-form we just have to consider $K_{n,q-1}$, the component of the kernel of bidegree $(0,q-1)$ in $z$ and $(n,n-q)$ in $\zeta$, that is:
\begin{align*}
\lefteqn{K_{n,q-1}(\zeta,z)}\\
&=\sum_{k=0}^{n-q} c_{n,k,q} \left(\frac{r(\zeta)}{r(\zeta)+S(\zeta,z)}\right)^{N+k}
\frac{\tilde s\wedge(\overline\partial_\zeta\tilde s)^{n-k-q}\wedge(\overline\partial_z\tilde s)^{q-1}\wedge(\overline\partial_\zeta\tilde Q)^k}{\langle \tilde s(\zeta,z),\zeta-z\rangle^{n-k}}.
\end{align*}
Moreover since the kernel is bounded when $|z-\zeta|$ is bounded away from zero and when $\zeta$ is away from the $bD$, we may assume that the form $f$ is supported in a neighborhood of $bD$ and that $z$ and $\zeta$ are close.\\
Let us fix $\zeta\in D$ near the boundary and set for $ j=1,2$
$$I_j(\zeta):=\int_{z\in \Omega_j} \left| f(\zeta)\wedge K_{n,q-1}(\zeta,z)\left[\frac{u_1(z)}{k(z, u_1(z))},\ldots,\frac{u_{q-1}(z)}{k(z, u_{q-1}(z))}\right] \right|d\lambda(z),$$
where $\Omega_1=\p{kd(\zeta)}\zeta\cap D$, $\Omega_2=D\setminus \p{kd(\zeta)}\zeta$ and $k>0$ is such that $\p{kd(\zeta)}\zeta\subset D_{\frac{r(\zeta)}2}$. We prove that $$I_j(\zeta)\leqs \sum_{i_1,\ldots, i_q}\left| d(\zeta)^{q-\sum_{i=n-q+2}^n\frac1{m_i}+\gamma }f(\zeta)\left[\frac{e_{i_1}(\zeta)}{k(\zeta, e_{i_1}(\zeta))},\ldots,\frac{e_{i_q}(\zeta)}{k(\zeta, e_{i_q}(\zeta))}\right]\right|$$
and then integrating with respect to $\zeta$ we get the estimate we are looking for.\\
For $I\subset\{1,\ldots, n\}$ we denote by $\overline{e_I}(\zeta)$ the family $\overline{e_i}(\zeta)$, $i\in I$ and by $\widehat{\overline{e_I}}(\zeta)$ the family $e_1(\zeta),\ldots, e_n(\zeta), \overline{e_1}(\zeta),\ldots, \overline{e_n}(\zeta)$ where the indices $i\in I$ have been removed. We also set $k(\zeta,\overline{ e_I}(\zeta))=\prod_{i\in I} k(\zeta,\overline{e_i}(\zeta))$.\\
We will work with forms which contain $d\zeta$ and $dz$ and we will compute them for some vectors $v$ and $w$. In order to indicate if the form which applies to $v$ is the form $dz$ or $d\zeta$ we write respectively $[v]_z$ or $[v]_\zeta$.\\
The two following lemmas give us estimates that we will need to evaluate $K(\zeta,z)[u_1,\ldots, u_{q-1}]_z [\widehat{\overline{e_I}}(\zeta)]_\zeta$.
\begin{lemma}\mlabel{lem1}
For all $\zeta\in D$ near enough $bD$, $\varepsilon\geq kd(\zeta)$, $z\in \p\varepsilon\zeta$, $w''_1,\ldots, w''_n$ a Yu basis at $z$ and $u=\sum_{j=1}^nu''_jw''_j\neq 0$ a vector of $\cc^n$, we have uniformly \wrt $\zeta,\varepsilon$ and $z$
$$\frac{\left| u''_j \right|}{\tau(\zeta,w''_j,\varepsilon)k(z,u)}\leqs \left(\frac{d(z)}{\varepsilon}\right)^{\frac{1}{m_j}} \frac{1}{d(z)}.$$
\end{lemma}
\pr By definition of $k(z,u)$ we have $\frac{\left| u''_j \right|}{\tau(\zeta,w''_j,\varepsilon)k(z,u)}=\frac{\left| u''_j \right|\tau(z,u,d(z))}{\tau(\zeta,w''_j,\varepsilon)d(z)}.$\\
On the one hand since $z\in\p\varepsilon\zeta$, by Proposition \ref{prop5} $\tau(z,w''_j,\varepsilon)\eqs\tau(\zeta,w''_j,\varepsilon)$.\\
On the other hand Proposition \ref{prop4} implies that $\frac{\left| u''_j \right|}{\tau(z,w''_j,\varepsilon)}\leqs \frac{1}{\tau(z,u,\varepsilon)}$. Thus Proposition \ref{prop4} yields  to 
$$\frac{\left| u''_j \right|}{\tau(\zeta,w''_j,\varepsilon)k(z,u)}\leqs\frac{\tau(z,w''_j,d(z))}{d(z)\tau(z,w''_j,\varepsilon)}\leqs \frac{1}{d(z)}\left( \frac{d(z)}{\varepsilon} \right)^{\frac{1}{m_j}}.\qed$$
\begin{lemma}\mlabel{lem2}
Let $\zeta\in D$ near enough $bD$, $q=1,\ldots, n$, $\varepsilon\geq kd(\zeta)$, $K\subset\{1,\ldots, n\}$, $I\subset\{1,\ldots, n\}$ such that ${\rm card}\ K=n-q$ and ${\rm card }\ I= q$.
Then we have uniformly \wrt $\zeta$ 
$$\left| \bigwedge_{i=1}^n d\zeta'_i \wedge \bigwedge_{k\in K} d\overline{\zeta'_k} \left(\widehat{\overline{e_I}}(\zeta)\right) k(\zeta,\overline{e_I}(\zeta))\right|\leqs \frac{\varepsilon^q Vol\bigl(\p{d(\zeta)}\zeta\bigr)}{\prod_{k\notin K} \tau'_k(\zeta,\varepsilon)}.$$
\end{lemma}
\pr The inequality   $\left| e'_{ij}(\zeta) \right|\leqs\tau'_j(\zeta,d(\zeta))$  given by Proposition \ref{prop10} yields to
\begin{align}
\left| \bigwedge_{i=1}^n d\zeta'_i \wedge \bigwedge_{k\in K} d\overline{\zeta'_k} \left(\widehat{\overline{e_I}}(\zeta)\right)\right|\leqs \prod_{i=1}^n\tau'_i(\zeta,d(\zeta)) \prod_{{k\in K}}\tau'_k(\zeta,d(\zeta))\mlabel{eq2}
\end{align}
Since $e_i(\zeta )$ is a unit vector for the Bergman metric, Proposition \ref{prop4} implies that $k(\zeta ,e_i(\zeta ))\eqs d(\zeta )$ 
which together with (\ref{eq2}) gives
\begin{align}
\left| \bigwedge_{i=1}^n d\zeta'_i \wedge \bigwedge_{k\in K} d\overline{\zeta'_k} \left(\widehat{\overline{e_I}}(\zeta)\right) k(\zeta,\overline{e_I}(\zeta))\right|\leqs \frac{d(\zeta)^q Vol\bigl(\p{d(\zeta)}\zeta\bigr)}{\prod_{k\notin K} \tau'_k(\zeta,d(\zeta))}.\mlabel{eq3}
\end{align}
Now, since $\varepsilon\geq kd(\zeta)$ we have $d(\zeta)\tau'_i(\zeta,d(\zeta))^{-1}\leqs 
\varepsilon\tau'_i(\zeta,\varepsilon)^{-1}$ for all $i$ which with (\ref{eq3}) proves the lemma.\qed\\
We now come to the estimate of $I_1$. Let $z$ be a point in $\p{kd(\zeta)}\zeta$. We denote by $u$ the family $u_1,\ldots, u_q$ and set $k(z,u(z))=\prod_{i=1}^{q-1} k(z, u_i(z))$, $\beta=q-1-\sum_{i=n-q+2}^n\frac{1}{m_i} +\gamma $ where $\gamma >-1$. We fix a Yu basis $w'_1,\ldots, w'_n$ at $\zeta$ and a Yu basis $w''_1,\ldots, w''_n$ at $z$ and we denote by $\zeta'$ and $z''$ the corresponding coordinates in the coordinates system centered at $\zeta$ and $z$.\\
Since $f\wedge K_{n,q-1}$ is of bidegree $(n,n)$ in $\zeta$ we have
\begin{align*}
\lefteqn{d(z)^\beta\left| f(\zeta)\wedge K_{n,q-1}(\zeta,z)\left[\frac{u(z)}{k(z,u(z))}\right]_z \right|=}\\
&=d(z)^\beta\left|f(\zeta)\wedge K_{n,q-1}(\zeta,z)\left[\frac{u(z)}{k(z,u(z))}\right]_z[e(\zeta),\overline{e}(\zeta)]_\zeta\right| \frac{1}{|\det(e(\zeta))|^2}\\
&\leq \sum_{\over{I\subset\{1,\ldots, n\}}{{\rm card } I= q}}  d(\zeta)^{\beta+1} \left|f(\zeta)\left[ \frac{\overline{e_I}(\zeta)}{k(\zeta,\overline{e_I}(\zeta))}\hskip -1pt \right]\right|\\
&  \hskip 30pt \cdot \frac{d(z)^\beta}{d(\zeta)^{\beta+1}}\left| K_{n,q-1}(\zeta,z) \left[\frac{u(z)}{k(z,u(z))} \right]_z [\widehat{\overline{e_I}}(\zeta)]_\zeta\right|  \frac{k(\zeta,\overline{e_I}(\zeta))}{|\det(e(\zeta))|^2}.
\end{align*}
Denote by ${\cal S}_n$ the set of all permutations of $\{1,\ldots, n\}$: $K_{n,q-1}$ is then a sum for $k=0,\ldots, n-q$ and $\nu,\mu,\lambda\in {\cal S}_n$ of
\begin{align*}
\left( \frac{r(\zeta)}{r(\zeta)+S(\zeta,z)} \right)^{N+k} \frac{\tilde K_{\nu,\mu,\lambda}(\zeta,z)}{\langle \tilde s(\zeta,z),\zeta-z\rangle^{n-k}}  \bigwedge_{i=1}^n d\zeta'_i\wedge\bigwedge_{i=q+1}^n d\overline{\zeta}'_{\mu(i)}\wedge\bigwedge_{i=2}^q d\overline z''_{\lambda(i)}
\end{align*}
where 
$$\tilde K_{\nu,\mu,\lambda}(\zeta,z)\hskip-1pt=\hskip-1pt \tilde s_{\nu(1)}(\zeta,z)\hskip-1pt \prod_{i=2}^q\diffp{\tilde s'_{\nu(i)}}{\overline{z}''_{\lambda(i)}}(\zeta,z)\hskip-1pt \prod_{i=q+1}^{n-k} \hskip-1pt\diffp{\tilde s'_{\nu(i)}}{\overline{\zeta}'_{\mu(i)}}(\zeta,z)\hskip-1pt \prod_{i=n-k+1}^{n}\hskip-1pt \diffp{\tilde Q'_{\nu(i)}}{\overline{\zeta}'_{\mu(i)}}(\zeta,z).$$
From Corollary \ref{cor3} we get
$$\left| \tilde K_{\nu,\mu,\lambda}(\zeta,z) \right|\leqs
\frac{ d(\zeta)^{4(n-k)} \sum_{j=1}^n\frac{\left| z'_j \right|}{\tau'_j(\zeta,d(\zeta))}}
{\prod_{i=1}^n \tau'_i(\zeta,d(\zeta)) \prod_{i=2}^q\tau(\zeta,w''_{\lambda(i)},d(\zeta))\prod_{i=q+1}^n\tau'_{\mu(i)}(\zeta,d(\zeta))}.$$
By Proposition \ref{prop11} $\left| \frac{r(\zeta)}{r(\zeta)+S(\zeta,z)} \right|\leqs 1$ and by Lemma \ref{lem3} $\left| \langle\tilde s(\zeta,z)\zeta-z\rangle \right|^{\frac12}\geqs d(\zeta)^2 \sum_{j=1}^n\frac{\left| z'_j \right|}{\tau'_j(\zeta,d(\zeta)}$.\\
Thus, keeping in mind that $d(\zeta)\eqs d(z)$,  with Lemmas \ref{lem1} and \ref{lem2} we get
\begin{align*}
\lefteqn{\frac{d(z)^\beta}{d(\zeta)^{\beta+1}} \left| K_{n,q-1}(\zeta,z)\left[\frac{u(z)}{k(z,u(z))}\right]_z \left[\widehat{\overline{e_I}}(\zeta) \right]_\zeta \right| \frac{k(\zeta,d(\zeta))}{\left| \det e(\zeta) \right|^2}}\\
&\leqs \sum_{k=0}^{n-q}\sum_{\lambda,\mu,\nu\in {\cal S}_n}\frac1{d(\zeta)}\frac{d(\zeta)^q Vol\bigl(\p{d(\zeta)}\zeta\bigr) }{\prod_{i=1}^q\tau'_{\mu(i)}(\zeta,d(\zeta))} \frac{1}{d(z)^{q-1}}\prod_{i=2}^q \left( \frac{d(z)}{d(\zeta)} \right)^{\frac{1}{m_{\lambda(i)}}}  \\
&\hskip 18pt \cdot\frac{ d(\zeta)^{4(n-k)} \sum_{j=1}^n\frac{\left| z'_j \right|}{\tau'_j(\zeta,d(\zeta))}}
{\prod_{i=1}^n \tau'_i(\zeta,d(\zeta)) \prod_{i=q+1}^n\tau'_{\mu(i)}(\zeta,d(\zeta))}
 \frac{Vol\bigl(\p{d(\zeta)}\zeta\bigr)^{-1}}{d(\zeta)^{4(n-k)} \left( \sum_{j=1}^n \frac{\left| z'_j \right|}{\tau'_j(\zeta,d(\zeta))} \right)^{2n-2k}}\\
&\leqs \frac{1}{\prod_{i=1}^n \tau'_i(\zeta,d(\zeta))^2\left( \sum_{j=1}^n\frac{\left| z'_j \right|}{\tau'_j(\zeta,d(\zeta))} \right)^{2n-1}}.
\end{align*}
Let $\alpha>0$ sufficiently small, which will be chosen in  moment. For $i=1,\ldots, n$ we have
$\tau'_i(\zeta,d(\zeta)) \sum_{j=1}^n \frac{\left| z'_j \right|}{\tau'_j(\zeta,d(\zeta))}\geqs \left| z'_i \right|$ and $d(\zeta)\eqs\tau'_1(\zeta,d(\zeta))$ so
\begin{align*}
\lefteqn{\frac{d(z)^\beta}{d(\zeta)^{\beta+1}} \left| K_{n,q-1}(\zeta,z)\left[\frac{u(z)}{k(z,u(z))}\right]_z \left[\widehat{\overline{e_I}}(\zeta) \right]_\zeta \right| \frac{k(\zeta,d(\zeta))}{\left| \det e(\zeta) \right|^2}}\\
&\hskip150pt\leqs\frac{d(\zeta)^{2(n-1)\alpha-1}}{ \left| z'_1 \right|^{1+2(n-1)\alpha}}\prod_{i=2}^{n} \frac{1}{\left| z'_i \right|^{2(1-\alpha)}\tau'_i(\zeta,d(\zeta))^{2\alpha}} .
\end{align*}
Now, if $2(n-1)\alpha+1<2$, integrating for $z\in\p{kd(\zeta )}\zeta $ we get
$$\int_{z\in\p{kd(\zeta)}\zeta} {\frac{d(z)^\beta}{d(\zeta)^{\beta+1}} \left| K_{n,q-1}(\zeta,z)\left[\frac{u(z)}{k(z,u(z))}\right]_z \left[\widehat{\overline{e_I}}(\zeta) \right]_\zeta \right| \frac{k(\zeta,d(\zeta))}{\left| \det e(\zeta) \right|^2}}\leqs 1$$
and so
$$I_1(\zeta )\leqs \sum_{\over{I\subset\{1,\ldots, n\}}{{\rm card } I=q}} d(\zeta )^{q-\sum_{i=n-q+2}^n\frac{1}{m_i} +\gamma } \left| f(\zeta )\left [\frac{\overline{e_I}(\zeta )}{k\bigl(\zeta ,\overline{e_I}(\zeta )\bigr)} \right]
 \right|.$$
Now we consider $I_2(\zeta )$. We cover the integration domain with the polyannuli $\po{2^jd(\zeta )}\zeta $, $j=0,1,\ldots$ and we consider for $\varepsilon \geq kd(\zeta )$ $$I_{2,\varepsilon }(\zeta ):=\int_{\po\varepsilon \zeta } \left| d(\zeta )^\beta  f(\zeta )\wedge K_{n,q-1}(\zeta ,z)\left [\frac{u(z)}{k\bigl(z,u(z)\bigr)} \right]_z \right|dV(z).$$
We proceed then in the same way as we did for $I_1$.\\
By Proposition \ref{prop11}, $\left| \frac{r(\zeta)}{r(\zeta)+S(\zeta,z)} \right|^{N+k}\leqs \left( \frac{d(\zeta)}{\varepsilon} \right)^{N+k}$ and from Corollary \ref{cor5} we get that $\left| \langle \tilde s(\zeta,z),\zeta-z\rangle \right|\geqs \varepsilon^4$. From Corollary \ref{cor3} we deduced that
\begin{align*}
\left| \tilde K_{\nu,\mu,\lambda}(\zeta,z) \right|
\hskip-1pt\leqs\hskip-1pt \frac{\varepsilon^{4(n-k)}}{\prod_{i=1}^n\tau'_i(\zeta,\varepsilon)
\prod_{i=2}^q\tau(\zeta,w''_{\lambda(i)},\varepsilon)\prod_{i=q+1}^n\tau'_{\mu(i)}(\zeta,\varepsilon)}\hskip-2pt\left( \frac{\varepsilon}{d(\zeta)} \right)^{2k}\hskip-3pt .
\end{align*}
Therefore Lemmas \ref{lem1} and \ref{lem2} put together yield to
\begin{align}
\nonumber\lefteqn{\frac{d(z)^\beta}{d(\zeta)^{\beta+1}} \left| K_{n,q-1}(\zeta,z)\left[\frac{u(z)}{k(z,u(z))}\right]_z \left[\widehat{\overline{e_I}}(\zeta) \right]_\zeta \right| \frac{k(\zeta,d(\zeta))}{\left| \det e(\zeta) \right|^2}}\\
\nonumber&\leqs 
\sum_{k=0}^{n-q}\sum_{\lambda,\mu,\nu\in {\cal S}_n}
\frac{d(z)^\beta}{d(\zeta)^{\beta+1}}
\frac{\varepsilon^q Vol\bigl(\p{d(\zeta)}\zeta\bigr) }{\prod_{i=1}^q\tau'_{\mu(i)}(\zeta,\varepsilon)} \frac{1}{d(z)^{q-1}}\prod_{i=2}^q \left( \frac{d(z)}{d(\zeta)} \right)^{\frac{1}{m_{\lambda(i)}}} \\
\nonumber&\hskip 18pt \cdot\frac{ \varepsilon^{4(n-k)} }
{\prod_{i=1}^n \tau'_i(\zeta,\varepsilon) \prod_{i=q+1}^n\tau'_{\mu(i)}(\zeta,\varepsilon))}
 \frac{Vol\bigl(\p{d(\zeta)}\zeta\bigr)^{-1}} {\varepsilon^{4(n-k)}}\left( \frac{d(\zeta)}{\varepsilon} \right)^{N-k}\\
\mlabel{eq71} &\leqs
\left( \frac{d(z)}{d(\zeta)} \right)^{\sum_{i=2}^q\frac1{m_{\lambda(i)}}-\sum_{i=n-q+2}^n\frac{1}{m_i}}
 \left( \frac{d(\zeta)}{\varepsilon} \right)^{N-n-\gamma} \frac{d(z)^\gamma \varepsilon^{-\gamma } }{Vol\ \p\varepsilon\zeta}.
\end{align}
Now $\sum_{i=2}^q \frac{1}{m_{\lambda(i)}}-\sum_{i=n-q+2}^n \frac{1}{m_i} \geq 0$ and since $z$ belongs to $\p\varepsilon\zeta$, we have $d(z)\leqs \varepsilon$ and so $\left( \frac{d(z)}{d(\zeta)} \right)^{\sum_{i=2}^q \frac{1}{m_{\lambda(i)}}-\sum_{i=n-q+2}^n \frac{1}{m_i}}\leqs \left( \frac{\varepsilon}{d(\zeta)} \right)^{q-1}$. Therefore
\begin{align*}
\lefteqn{\frac{d(z)^\beta}{d(\zeta)^{\beta+1}} \left| K_{n,q-1}(\zeta,z)\left[\frac{u(z)}{k(z,u(z))}\right]_z \left[\widehat{\overline{e_I}}(\zeta) \right]_\zeta \right| \frac{k(\zeta,d(\zeta))}{\left| \det e(\zeta) \right|^2}}\\
&\hskip 130pt\hbox{\hfill} \leqs d(z)^\gamma \varepsilon^{-\gamma } \left( \frac{d(\zeta)}{\varepsilon} \right)^{N-n-q+1-\gamma} \frac{1}{Vol\ \p\varepsilon\zeta}
\end{align*}
and so when we integrate with respect to $z\in \po\varepsilon\zeta$: $$I_{2,\varepsilon}(\zeta)\leqs \left( \frac{d(\zeta)}{\varepsilon} \right)^{N-n-q+1-\gamma} \sum_{\over{I\subset\{1,\ldots, n\}}{{\rm card } I=q}} d(\zeta )^{\beta+1} \left| f(\zeta )\left [\frac{\overline{e_I}(\zeta )}{k\bigl(\zeta ,\overline{e_I}(\zeta )\bigr)} \right]
 \right|.$$
Provided $N$ is large enough we conclude that 
\begin{align*}
I_2(\zeta)&\leq \sum_{j=0}^{+\infty}I_{2,2^jkd(\zeta)}(\zeta)\\
&\leqs\sum_{j=0}^{+\infty} \left( \frac{1}{2^jkd(\zeta)} \right)^{N-n-q+1-\gamma} \hskip -9pt  \sum_{\over{I\subset\{1,\ldots, n\}}{{\rm card } I=q}} d(\zeta )^{\beta+1} \left| f(\zeta )\left [\frac{\overline{e_I}(\zeta )}{k\bigl(\zeta ,\overline{e_I}(\zeta )\bigr)} \right]
 \right|\\
&\leqs \sum_{\over{I\subset\{1,\ldots, n\}}{{\rm card } I=q}} d(\zeta )^{\beta+1} \left| f(\zeta )\left [\frac{\overline{e_I}(\zeta )}{k\bigl(\zeta ,\overline{e_I}(\zeta )\bigr)} \right]
 \right|
\end{align*}
and finaly that Theorem \ref{th1} holds for $\alpha=0$, which proves the $W^0$ estimates.
We will use the same kind of method for the $W^1$-estimates as for the $W^0$-estimates and we keep the same notation. We have to show that for every $z_0\in bD$ and every sufficiently small $\varepsilon_0>0$
\begin{align}
\int_{z\in\p{\varepsilon_0}{z_0}\cap bD} \left| d(z)^\beta Tf(z) \right| dV(z)\leqs \sigma(\p{\varepsilon_0}{z_0}\cap bD) \left\|d(\cdot)^{\beta+1} f(\cdot)\right\|_{\w1q}\mlabel{eq4} 
\end{align}
where as before $\beta=q-1+\sum_{i=n-q+2}^n \frac{1}{m_i}+\gamma,$ $\gamma\in]-1,\gamma_0]$.\\
Since $\delta$ is a pseudodistance, there exists $K>0$ so large that $\delta(\zeta,z_0)\eqs \delta(\zeta,z)$ for all $\zeta\notin \p{K\varepsilon_0}{z_0}$ and all $z\in \p{\varepsilon_0}{z_0}$, $K$ not depending on $z_0$, $\zeta$, $z$ nor on $\varepsilon_0 $. We set for $j=1,2$ $$J_i(\zeta)=\int_{z\in \tilde \Omega_j }d(z)^\beta\left| f(\zeta)\wedge K_{n,q-1}(\zeta,z)
\left[\frac{u(z)}{k\bigl(z,u(z)\bigr)} \right]_z
 \right|dV(z)$$
where $\tilde \Omega_1= \p{kd(\zeta)}\zeta\cap \p{\varepsilon_0}{z_0}\cap D$ and
$\tilde \Omega_2=  \bigl(\p{\varepsilon_0}{z_0}\cap D\setminus \p{kd(\zeta)}\zeta\bigr)$. We estimate both integrals.\\
From the $W^0$-estimates we have for $i=1,2$
$$J_i(\zeta)\leqs \sum_{\over{I\subset\{1,\ldots, n\}}{{\rm card } I=q}} d(\zeta )^{\beta+1} \left| f(\zeta )\left [\frac{\overline{e_I}(\zeta )}{k\bigl(\zeta ,\overline{e_I}(\zeta )\bigr)} \right]
 \right|$$
and therefore, since $\sigma(\p{K\varepsilon_0}{z_0}\cap bD) \eqs \sigma(\p{\varepsilon_0}{z_0}\cap bD)$ we obtain that
\begin{align*}
\lefteqn{\int_{\zeta\in\p{K\varepsilon_0}{z_0}}
\int_{z\in \tilde \Omega_j }d(z)^\beta\left| f(\zeta)\wedge K_{n,q-1}(\zeta,z)
\left[\frac{u(z)}{k\bigl(z,u(z)\bigr)} \right]_z
 \right|dV(z) dV(\zeta)}\\
&\hskip 66pt \leqs\sum_{\over{I\subset\{1,\ldots, n\}}{{\rm card } I=q}} 
\int_{\zeta\in\p{K\varepsilon_0}{z_0}}
d(\zeta )^{\beta+1} \left| f(\zeta )\left [\frac{\overline{e_I}(\zeta )}{k\bigl(\zeta ,\overline{e_I}(\zeta )\bigr)} \right]
 \right| dV(\zeta)\\
&\hskip 66pt \leqs 
\sigma(\p{K\varepsilon_0}{z_0}\cap bD) \left\|d(\cdot)^{\beta+1} f(\cdot)\right\|_{\w1q}\\
&\hskip 66pt \eqs \sigma(\p{\varepsilon_0}{z_0}\cap bD) \left\|d(\cdot)^{\beta+1} f(\cdot)\right\|_{\w1q}.
\end{align*}
Now let us consider the case where $\zeta\in D\cap\po\varepsilon {z_0}$, $\varepsilon\geq K\varepsilon_0$.\\
We have $\delta(\zeta,z)\eqs\delta(\zeta,z_0) \eqs \varepsilon$ so from Corollary \ref{cor5}, Proposition \ref{prop11} and Remark \ref{rem1} we have $\left| \frac{r(\zeta)}{r(\zeta)+S\left( \zeta,z \right)} \right|\leqs \frac{d(\zeta)}\varepsilon$ and $\left| \langle \tilde s(\zeta,z),\zeta-z\rangle \right|\geqs \varepsilon^4$. As in the case of the $W^0$-estimates we get 
\begin{align*}
\lefteqn{{d(z)^\beta} \left|f(\zeta)\wedge K_{n,q-1}(\zeta,z)\left[\frac{u(z)}{k(z,u(z))}\right]_z  \right|}\\
&\leqs 
\sum_{\over{I\subset\{1,\ldots, n\}}{{\rm card } I=q}} d(\zeta )^{\beta+1} \left| f(\zeta )\left [\frac{\overline{e_I}(\zeta )}{k\bigl(\zeta ,\overline{e_I}(\zeta )\bigr)} \right]
 \right|\\
& \cdot\sum_{\lambda\in {\cal S}_n}\hskip -5pt
\left( \frac{d(z)}{d(\zeta)} \right)^{
\sum_{i=2}^q \frac{1}{m_{\lambda(i)}}-\sum_{i=n-q+2}^n\frac{1}{m_i}}
\hskip -2pt\left( \frac{d(\zeta)}{\varepsilon} \right)^{N-n-\gamma}\hskip -3pt
 \frac{d(z)^\gamma  \varepsilon^{-\gamma }}{Vol\ \p\varepsilon\zeta}.
\end{align*}
Since $\zeta$ belongs to  $\p\varepsilon{z_0}$ Proposition \ref{prop1.2} gives   that $Vol\ \p\varepsilon\zeta \eqs Vol\ \p\varepsilon{z_0}$. Morevover $\varepsilon\geq K\varepsilon_0\geqs d(z)$ because $z$ belongs to $\p{\varepsilon_0}{z_0}$, so provided $N$ is large enough
\begin{align*}
\lefteqn{{d(z)^\beta} \left|f(\zeta)\wedge K_{n,q-1}(\zeta,z)\left[\frac{u(z)}{k(z,u(z))}\right]_z  \right|}\\
&\leqs \sum_{\over{I\subset\{1,\ldots, n\}}{{\rm card } I=q}} d(\zeta )^{\beta+1} \left| f(\zeta )\left [\frac{\overline{e_I}(\zeta )}{k\bigl(\zeta ,\overline{e_I}(\zeta )\bigr)} \right]
 \right|
d(z)^\gamma  \varepsilon^{-\gamma } \frac{1}{Vol\bigl(\p\varepsilon{z_0} \bigr)}.
\end{align*}
We integrate \wrt $z$ and get
$$J_2(\zeta)\leqs \frac{\sigma\bigl(\p{\varepsilon_0}{z_0}\cap bD\bigr)}{Vol\bigl(\p\varepsilon{z_0}\bigr)} \frac{\varepsilon_0^{\gamma +1}}{\varepsilon^\gamma }\sum_{\over{I\subset\{1,\ldots, n\}}{{\rm card } I=q}} d(\zeta )^{\beta+1} \left| f(\zeta )\left [\frac{\overline{e_I}(\zeta )}{k\bigl(\zeta ,\overline{e_I}(\zeta )\bigr)} \right]
 \right|$$
so 
$$\int_{\zeta\in\po\varepsilon{z_0}\cap D} J_2(\zeta) dV(\zeta) 
\leqs {\sigma\bigl(\p{\varepsilon_0}{z_0}\cap bD\bigr)} \left( \frac{\varepsilon_0}{\varepsilon} \right)^{\gamma +1}\left\| d(\cdot)^{\beta+1} f(\cdot)\right\|_{\w1q}.$$
Now since $\gamma >-1$:
\begin{align*}
\int_{\zeta\in D\setminus \p{K\varepsilon_0}{z_0} } J_2(\zeta)dV(\zeta)&\leqs
\sum_{j=0}^{+\infty}
\int_{\zeta\in D\cap \po{2^jK\varepsilon_0}{z_0} } J_2(\zeta)dV(\zeta)\\
&\leqs {\sigma\bigl(\p{\varepsilon_0}{z_0}\cap bD\bigr)}  \left\| d(\cdot)^{\beta+1} f(\cdot)\right\|_{\w1q}.
\end{align*}
Let us notice that if $k$ is small enough $\tilde \Omega_1=\p{kd(\zeta)}\zeta\cap bD\cap \p{\varepsilon_0}{z_0}=\emptyset$ for all $\zeta\notin \p{K\varepsilon_0}{z_0}$. Indeed, for $z\in \tilde \Omega_1$ we have $d(z)\eqs d(\zeta )$ and $\delta (\zeta ,z)\leq kd(\zeta )$ because $z$ belongs to $\p{kd(\zeta )}\zeta$ and $d(z)\leqs \varepsilon _0$ because $z$ also belongs to  $\p{\varepsilon _0}{z_0}$ so $k\varepsilon _0\geqs \delta (\zeta ,z)$. Morevover $\delta (\zeta ,z)\eqs \delta (\zeta ,z_0)\geq K\varepsilon _0$ for all $\zeta\in\p{\varepsilon _0}{z_0}$ so $k\varepsilon _0\geqs K\varepsilon _0$ and it is impossible if $k$ is small enough. Therefore $J_1(\zeta)=0$ for all  $\zeta \notin \p{K\varepsilon _0}{z_0}$ and we thus have proved that (\ref{eq4}) holds true. Therefore Theorem \ref{th3} holds for $\alpha=0$ and $\alpha=1$ so, by interpolation it  holds for all $\alpha\in[0,1]$ which finishes the proof of Theorem \ref{th3}.
\section{The limit case}
We are now ready to prove Theorem \ref{th2}. At first we consider the BMO-case and the $L^1$-case. Again we just have to consider the case of smooth forms. Let $f\in C^\infty_{0,1}(\overline D)$. We set $u(z)=\int_D f(\zeta)\wedge K_{n,0}(\zeta,z)$ and show that $u$ is continuous up to the boundary  by proving that for every $z\in bD$, $K(\cdot, z)$ is uniformly integrable (see \cite{Sko}).\\
For $z\in bD$ we have $\tilde s(\zeta,z)=s(\zeta,z).$ 
From Proposition \ref{prop6} we easily obtain:
\begin{proposition}\mlabel{prop20}
For all $\zeta\in D$ near $bD$, all $\varepsilon>0$ and all $z\in \p\varepsilon\zeta\cap bD$ we have uniformly \wrt $\zeta$, $z$ and $\varepsilon$
\begin{align*}
\left|  s'_j(\zeta,z) \right|&\leqs\frac{\varepsilon^4}{\tau'_j(\zeta,\varepsilon)},\\
\left| \diffp{s'_j}{{z'_i}}(\zeta,z) \right|,\left| \diffp{s'_j}{\overline{\zeta'_i}}(\zeta,z) \right|, \left| \diffp{ s'_j}{\overline{z'_i}}(\zeta,z) \right|&\leqs \frac{\varepsilon^4}{\tau'_i(\zeta,\varepsilon) \tau'_j(\zeta,\varepsilon)},\\
 \left| \diffp{ s'_j}{\overline{z'_i}\partial z'_k}(\zeta,z) \right|,\left| \diffp{s'_j}{\overline{z'_i} \partial \overline{z'_k}}(\zeta,z) \right|
&\leqs \frac{\varepsilon^4}{\tau'_i(\zeta,\varepsilon) \tau'_j(\zeta,\varepsilon)\tau'_k(\zeta,\varepsilon)},\\
\left| \diffp{s'_j}{\overline{\zeta'_i} \partial z'_k}(\zeta,z) \right|,
\left| \diffp{s'_j}{\overline{\zeta'_i} \partial \overline z'_k}(\zeta,z) \right|&\leqs \frac{\varepsilon^4}{\tau'_i(\zeta,\varepsilon) \tau'_j(\zeta,\varepsilon)\tau'_k(\zeta,\varepsilon)}.
\end{align*}
\end{proposition}
And so for all $\varepsilon\geq kd(\zeta)$ and all $z\in bD\cap \po\varepsilon\zeta$ we have as in the case of $W^0$-estimates
\begin{align*}
\left| K_{n,0}(\zeta,z) \right|&\leqs \sum_{j=1}^n \sum_{k=0}^{n-1}
\left( \frac{d(\zeta)}{\varepsilon} \right)^{N+k}
\frac{\varepsilon^{4(n-k)}}{\prod_{i=1}^n\tau'_i(\zeta,\varepsilon)^2 \varepsilon^{4(n-k)}} \left( \frac{\varepsilon}{d(\zeta)} \right)^{2k}\tau'_j(\zeta,\varepsilon)\\
&\leqs \frac{\varepsilon^{\frac{1}{m}}}{Vol\bigl(\p\varepsilon\zeta\bigr)}\eqs \frac{\varepsilon^{\frac{1}{m}}}{Vol\bigl(\p\varepsilon z\bigr)}.
\end{align*}
Therefore $\int_{D\cap \p\varepsilon z}\left| K_{n,0}(\zeta,z) \right|dV(\zeta)=O(\varepsilon^{\frac1m})$ and equation (24) of \cite{BCD} is satisfied. We conclude as in \cite{BCD} and \cite{DM} that $K_{n,0}(\cdot,z)$ is uniformly integrable and that $u$ is continuous up to $bD$.\\
Now to prove the $L^1$ estimates we aim to show that
\begin{align*}
\int_{z\in bD} \left| f(\zeta)\wedge K_{n,0}(\zeta,z) \right| d\sigma(z) \leqs \sum_{j=1}^n \left| f(\zeta)\left[\frac{\overline{e_j}(\zeta)}{k(\zeta,\overline{e_j}(\zeta)} \right] \right|.
\end{align*}
We have 
\begin{align*}
\left| f(\zeta)\wedge K_{n,0}(\zeta,z) \right|&\leqs
\sum_{j=1}^n \left| f(\zeta)\left[\frac{\overline{e_j}(\zeta)}{k(\zeta,\overline{e_j}(\zeta)} \right] \right| \left| K_{n,0}(\zeta,z)[\widehat{\overline{e_j}}(\zeta) ] \right|\frac{k\bigl(\zeta,\overline{e_j}(\zeta)\bigr)}{\left| \det e(\zeta) \right|^2}.
\end{align*}
From Propositions \ref{prop11}, \ref{prop10}, \ref{prop20}, Corollary \ref{cor5} and Lemma \ref{lem2} we get for all $z\in \po\varepsilon\zeta\cap bD$
\begin{align}
\left| K_{n,0}(\zeta,z)[\widehat{\overline{e_j}}(\zeta) ] \right|\frac{k\bigl(\zeta,\overline{e_j}(\zeta)\bigr)}{\left| \det e(\zeta) \right|^2}&\leqs \frac{\varepsilon}{Vol\bigl(\p\varepsilon\zeta\bigr)} \left( \frac{d(\zeta)}{\varepsilon} \right)^{N-n}\mlabel{eq6}
\end{align}
We therefore have
\begin{align*}
\int_{bD} 
\left| f(\zeta)\wedge K_{n,0}(\zeta,z) \right|d\sigma (z)&\leqs \sum_{j=1}^n \sum_{i=0}^{+\infty}
\int_{\po{2^i kd(\zeta)}\zeta\cap bD}\hskip -3pt
\left| f(\zeta)\wedge K_{n,0}(\zeta,z) \right|d\sigma (z)\\
&\leqs \sum_{j=1}^n \sum_{i=0}^{+\infty}\left( \frac{1}{2^jk} \right)^{N-n}
\left| f(\zeta)\left[\frac{\overline{e_j}(\zeta)}{k(\zeta,\overline{e_j}(\zeta)} \right] \right|\\
&\leqs \sum_{j=1}^n\left| f(\zeta)\left[\frac{\overline{e_j}(\zeta)}{k(\zeta,\overline{e_j}(\zeta)} \right] \right|
\end{align*}
and now integrating with respect to $\zeta$ we get $\int_{bD}|u(z)|d\sigma(z) \leqs \|f\|_{\w01}$ which proves the $L^1$-estimates.\\
Now for the BMO estimates we set for $j=1,2$
$$I_j\hskip -2pt=\hskip -4pt \int_{(z,w)\in (\p{\varepsilon_0}{z_0} \cap bD)^2}\hskip -3pt \int_{\zeta\in \tilde \Omega_j}\hskip -6pt \left|\hskip -1pt f(\zeta)\hskip -2pt\wedge\hskip -2pt \bigl(\hskip -1pt K_{n,0}(\zeta,z)\hskip -1pt-\hskip -1pt K_{n,0}(\zeta,w)\hskip -1pt\bigr)\hskip -1pt \right|\hskip -2pt dV\hskip -1pt(\zeta)d\sigma(z)d\sigma(w)$$
where $\tilde \Omega_1= D\cap \p{K\varepsilon_0}{z_0}$, $\tilde \Omega_2= D\setminus \p{K\varepsilon_0}{z_0}$, $K$ as for the $W^1$-estimates is such that $\delta(\zeta,z_0)\eqs\delta(\zeta,z)$ for all $z\in\p{\varepsilon_0}{z_0}$ and all $\zeta\notin\p{K\varepsilon_0}{z_0}$. 
To prove the BMO-estimates, it suffices to prove that, up to a multiplicative constant not depending on $z_0$ nor on $\varepsilon_0$, $I_1$ and $I_2$ are less or equal than  $\sigma(\p{\varepsilon_0}{z_0}\cap bD)^2 \|f\|_{\w11}$. We first consider $I_1$.\\
It suffices to show that
$\tilde I_1 =\int_{z\in\p{\varepsilon_0}{z_0}\cap bD} \int_{\tilde \Omega_1} \left| f(\zeta)\wedge K_{n,0}(\zeta,z) \right|dV(\zeta)d\sigma(z)$ is controlled by $\sigma(\p{\varepsilon_0}{z_0}\cap bD) \|f\|_{\w11}$.
From equation (\ref{eq6}) of the $L^1$ case we get that for all $\zeta\in\p{K\varepsilon_0}{z_0}$
\begin{align*}
\lefteqn{\int_{z\in\p{\varepsilon_0}{z_0}\cap bD} \left| f(\zeta)\wedge K_{n,0}(\zeta,z) \right|dV(\zeta)d\sigma(z)}\\
&\leqs \sum_{j=1}^n \left| f(\zeta)\left[ \frac{\overline e_j(\zeta)}{k(\zeta,\overline{e}_j(\zeta))}\right] \right|  \sum_{i=0}^\infty \int_{z\in\po{2^j k d(\zeta)}\zeta\cap bD} \frac{(2^i k)^{n-N} 2^ikd(\zeta)}{Vol \bigl(\p{2^jkd(\zeta)}\zeta\bigr)}d\sigma(z)\\
&\leqs \sum_{j=1}^n \left| f(\zeta)\left[ \frac{\overline e_j(\zeta)}{k(\zeta,\overline{e}_j(\zeta))}\right] \right|
\end{align*}
thus $I_1\leqs \sigma(\p{K\varepsilon_0}{z_0}\cap bD) \|f\|_{\w11}\eqs \sigma(\p{\varepsilon_0}{z_0}\cap bD)\|f\|_{\w01}$.\\
It is a little more difficult to handle $I_2$. We will evaluate $I_2$ in the way as we did for $J_2$ in the $W^1$-estimates. We fix $\varepsilon\geq K\varepsilon_0$ and set for $\zeta,z\in D$
$$\hat K_{n,0}(\zeta,z)= 
\sum_{k=0}^{n-1} c_{n,k} \left(\frac1{1+\tilde S(\zeta ,z)}\right)^{N+k} \frac{s \wedge (\overline\partial_\zeta  \tilde Q)^k\wedge (\overline \partial_\zeta  s)^{n-k-1}}{\langle  s(\zeta ,z),\zeta -z\rangle^{n-k}}.$$
We have $\hat K_{n,0}(\zeta,z)=K_{n,0}(\zeta,z)$ for all $z\in bD$ and for $w,z\in bD$, $z\neq w$:
$$K_{n,0}(\zeta,z)-K_{n,0}(\zeta,w)=\int_0^1 \diffp{\hat K_{n,0}}v\bigl(\zeta,w+t(z-w)\bigr) |z-w| dt$$
where $v=\frac{z-w}{\left| z-w \right|}$.\\
For $\zeta\in \po\varepsilon{z_0}$, $z,w\in\p{\varepsilon_0}{z_0}$ and $t\in[0,1]$ we have  $\delta(\zeta,w+t(z-w))\eqs\delta(\zeta,z_0)\eqs\varepsilon$. Proposition \ref{prop11} therefore holds and gives $|\langle s(\zeta,w+t(z-w)), w+t(z-w)\rangle |\eqs \varepsilon^4$ and $|r(\zeta)+S(\zeta,w+t(z-w))|\eqs \varepsilon$. Combining these inequalities with those of Propositions \ref{prop6} and \ref{prop4} we get as previously that for all
$z,w\in\p{\varepsilon_0}{z_0}$ and $t\in[0,1]$:
\begin{align}
\left|\hskip -1pt \diffp{\hat K_{n,0}}{v}(\zeta,w\hskip-1pt+\hskip-1pt t(z\hskip-1pt -\hskip-1pt w))[\widehat{\overline{e}_j}(\zeta)]k(\zeta,\overline{e}_j(\zeta)) \right|\hskip -2pt |z-w|\hskip-1pt
\leqs\hskip-1pt \frac{d(\zeta)^{N-n} \varepsilon^{1+n-N}}{Vol\ \p{\varepsilon}\zeta} \frac{|z-w|}{\tau(\zeta,v,\varepsilon)}.\mlabel{eq7}
\end{align}
On the one hand $\zeta$ belongs to $\p\varepsilon{z_0}$ so $\tau(\zeta,v,\varepsilon)\eqs\tau(z_0,v,\varepsilon)$. On the other hand $z,w$ belong to $\p{\varepsilon_0}{z_0}$ so $|z-w|\leqs \tau(z_0,v,\varepsilon_0)$ and since $\varepsilon_0\leq\varepsilon$ we have $\frac{\tau(z_0,v,\varepsilon_0)}{\tau(z_0,v,\varepsilon)}\leqs \left( \frac{\varepsilon_0}{\varepsilon} \right)^{\frac{1}{m}}$. Moreover, $
Vol\bigl(\p\varepsilon\zeta\bigr)\eqs Vol\bigl(\p\varepsilon{z_0}\bigr)$ and $d(\zeta)\leqs \varepsilon$ because $\zeta$ belongs to $\po{\varepsilon}{z_0}$. 
Inequality (\ref{eq7}) then becomes
\begin{align}
\left| \diffp{\hat K_{n,0}}{v}(\zeta,w+t(z-w))[\widehat{\overline{e}_j}(\zeta)]k(\zeta,\overline{e}_j(\zeta)) \right|
\leqs \left( \frac{\varepsilon_0}{\varepsilon} \right)^{\frac{1}{m}} \frac{d(\zeta)^{N-n} \varepsilon^{n+1-N}}{Vol \ \p{\varepsilon}{z_0}}
\end{align}
and integrating with respect to $z$ and $w$:
\begin{align*}
\lefteqn{\int_{z,w\in\p{\varepsilon_0}{z_0}\cap bD} 
\left|\bigl(K_{n,0}(\zeta,z)-K_{n,0}(\zeta,w)\bigr) [\widehat{\overline{e}_j}(\zeta)]k(\zeta,\overline{e}_j(\zeta))\right|}\\
&\hskip 150pt \leqs \left( \frac{\varepsilon_0}{\varepsilon} \right)^{\frac{1}{m}} \frac{1}{\sigma(\p\varepsilon{z_0}\cap bD)} {\sigma(\p{\varepsilon_0}{z_0}\cap bD)}^2.
\end{align*}
Integrating for $\zeta\in \po{2^iK\varepsilon_0}{z_0}$ and summing from $i=0$ to infinity we get
\begin{align*}
I_2\hskip -2pt&\leqs\hskip -2pt\sum_{j=1}^n\hskip -2pt \sum_{i=0}^\infty \hskip -2pt \int_{\zeta\in\po{2^iK\varepsilon_0}{z_0}}\hskip -2pt
\frac{ {2^{-\frac im}\sigma(\p{\varepsilon_0}{z_0}\cap bD)}^2} {\sigma(\p{2^{i} K \varepsilon_0}{z_0}\cap bD)} \left| f(\zeta)\hskip -1pt\left[\hskip -2pt\frac{\overline{e}_j(\zeta)}{k(\zeta,\overline{e}_j(\zeta))}\hskip -2pt\right]\hskip -1pt \right|dV(\zeta)\\
&\leqs\hskip -2pt {\sigma(\p{\varepsilon_0}{z_0}\cap bD)}^2\|f\|_{\w11}.
\end{align*}
This ends the proof of the BMO-estimates. Therefore  the limit case is proved for $L^1$ and BMO classes. By interpolation between $L^1$ and BMO we can conclude that it holds true for all $L^p$ classes with $p\in[1,+\infty[$. However we were unable to locate in the litterature a proper reference for interpolation results between $L^1$ and BMO spaces over convex domains of finite type and the homogeneous space structure we use on these domains. Hence for completeness' sake we briefly indicate how to conclude without interpolation results between these spaces.
We set for $j\in \{1,\ldots, n\},$ $\zeta,z\in bD$ and $t>0$
\begin{align*}
 \tilde K_t(\zeta ,z)&:=K_{n,0}(\zeta -t\eta_\zeta ,z) \cdot \frac{k\bigl(\zeta ,\overline{e_j}(\zeta )\bigr)}{\left| \det (e(\zeta )\right|^2},\\
\tilde K^0_t(\zeta ,z)&:= \frac1{\sigma (\p t \zeta \cap bD)} \chi_{\p t \zeta \cap bD} (z)
\end{align*}
where $\chi_{\p t \zeta \cap bD}$ is the characteristic function of $\p t \zeta \cap bD$. We also define for $\xi\in bD$ the following maximal functions:
\begin{align*}
Mf(\xi)&:= \sup_{\over{\zeta \in bD, t>0}{\delta(\zeta ,\xi)<t}} \int_{z\in \p t\zeta } f(z)\tilde K_t(\zeta ,z)d\sigma (z),\\
M^0f(\xi)&:= \sup_{\over{\zeta \in bD, t>0}{\delta(\zeta ,\xi)<t}} \int_{z\in \p t\zeta } f(z)\tilde K^0_t(\zeta ,z)d\sigma (z).
\end{align*}
$M^0f$ is the  Hardy-Littlewood maximal function associated to the homogeneous space structure of $bD$. From inequality (\ref{eq6}), which is the analog of Hypothesis (H2) of \cite{AB}, one can easily show that $Mf(\xi)\leqs M^0f(\xi)$ uniformaly with respect to $\xi$ and $f$ (see Remark 1, Section 4 of \cite{AB}). Thus since $\tilde K^0_t$ satisfies Hypothesis (H1) of \cite{AB}, $\tilde K_t$ also satisfies Hypothesis (H1). Theorems 2 and 3 from \cite{AB} imply that for all $f\in \w\alpha1(D)$, all $j=1,\ldots, n$
$\int_{\zeta\in bD}  
\left| K_{n,0}(\zeta ,z) [\widehat{\overline{e}_j}(\zeta)]k(\zeta,\overline{e}_j(\zeta)) \right| 
\left|f(\zeta )\left[\frac{\overline{e_j}(\zeta )}{k(\zeta , \overline{e_j}(\zeta )}\right]\right|dV(\zeta )$ belongs to $L^p(bD)$, $p=\frac1{1-\alpha},$
and has norm controled by $\| f\|_{\w\alpha1}$. Hence the limit case holds true for all $p\in]1,+\infty[$ which finishes the proof of theorem \ref{th2}.

\end{document}